# On the R-superlinear convergence of the KKT residues generated by the augmented Lagrangian method for convex composite conic programming

**Ying Cui · Defeng Sun · Kim-Chuan Toh**

**Abstract** Due to the possible lack of primal-dual-type error bounds, the superlinear convergence for the Karush-Kuhn-Tucker (KKT) residues of the sequence generated by augmented Lagrangian method (ALM) for solving convex composite conic programming (CCCP) has long been an outstanding open question. In this paper, we aim to resolve this issue by first conducting convergence rate analysis for the ALM with Rockafellar's stopping criteria under only a mild quadratic growth condition on the dual of CCCP. More importantly, by further assuming that the Robinson constraint qualification holds, we establish the R-superlinear convergence of the KKT residues of the iterative sequence under easy-to-implement stopping criteria for the augmented Lagrangian subproblems. Equipped with this discovery, we gain insightful interpretations on the impressive numerical performance of several recently developed semismooth Newton-CG based ALM solvers for solving linear and convex quadratic semidefinite programming.

The research of the second author was supported in part by the Academic Research Fund (grant R-146-000-207-112).

Ying Cui
Department of Mathematics, National University of Singapore, 10 Lower Kent Ridge Road, Singapore
E-mail: matcuiy@nus.edu.sg

Defeng Sun
Department of Mathematics and Risk Management Institute, National University of Singapore, 10 Lower Kent Ridge Road, Singapore
E-mail: matsundf@nus.edu.sg

Kim-Chuan Toh
Department of Mathematics and Institute of Operations Research and Analytics, National University of Singapore, 10 Lower Kent Ridge Road, Singapore
E-mail: mattohkc@nus.edu.sg

# 1 Introduction

Since the introduction of the augmented Lagrangian method (ALM) by Hestenes [28] and Powell [42] in the late 1960s for solving equality constrained optimization problems, the study on the ALM has grown into a fruitful subject in optimization, supported by significant theoretical developments over the past half-century; see, e.g., the papers [31,12,13,11,41,55,56] and the monographs [7,25]. The research of the ALM is advanced by its impressive numerical performance for various applications, including optimal control, partial differential equations and game theory [23,6,39,4]. Recently, the algorithm has also been successfully implemented in several efficient solvers for large scale convex positive semidefinite programming (SDP) [60,59,36] and convex composite matrix programming [33,10].

The research of the present paper is motivated by our desire to deeply understand the numerical success of the ALM in the aforementioned convex programming solvers. Notably, the targeted problems of these solvers belong to a wide class of convex programming – convex composite conic programming (CCCP), for which the objective function is the sum of a nonsmooth term and a smooth term, and the constraints include linear equations and convex cones. More specifically, a CCCP problem takes the form of

$$\begin{aligned} \min \ & f^0(x) := h(\mathcal{A}x) + \langle c, x \rangle + p(x) \\ \text{s.t.} \ & \mathcal{B}x \in b + \mathcal{Q}, \end{aligned} \tag{P}$$

where $\mathcal{A} : \mathbb{X} \to \mathbb{W}$ and $\mathcal{B} : \mathbb{X} \to \mathbb{Y}$ are linear maps, $\mathcal{Q} \subseteq \mathbb{Y}$ is a closed convex cone, $c \in \mathbb{X}$ and $b \in \mathbb{Y}$ are given data, $p : \mathbb{X} \to (-\infty, +\infty]$ is a proper closed convex function, $h : \mathbb{W} \to (\infty, +\infty]$ is a proper closed convex and essentially smooth function, whose gradient is locally Lipschitz continuous on $\text{int}(\text{dom}\, h)$, and $\mathbb{X}, \mathbb{Y}$ and $\mathbb{W}$ are three finite dimensional Euclidean spaces. In particular, when $\mathcal{Q}$ is the cone of symmetric and positive semidefinite matrices, (P) reduces to a convex SDP problem; when $p$ is the matrix norm function, (P) reduces to a convex composite matrix optimization problem.

Let $\sigma > 0$ be a given penalty parameter. The augmented Lagrangian function associated with (P) is given by (cf. [52, Section 11.K] or [55])

$$L_\sigma(x,y) := f^0(x) + \frac{1}{2\sigma} \big( \|\Pi_{\mathcal{Q}^\circ}[y + \sigma(\mathcal{B}x - b)]\|^2 - \|y\|^2 \big), \quad (x,y) \in \mathbb{X} \times \mathbb{Y},$$

where $\mathcal{Q}^\circ \subset \mathbb{Y}$ is the polar cone of $\mathcal{Q}$, i.e., $\mathcal{Q}^\circ = \{y \in \mathbb{Y} \,|\, \langle y, d \rangle \leqslant 0, \ \forall d \in \mathcal{Q}\}$, and $\Pi_{\mathcal{Q}^\circ}(\cdot)$ denotes the Euclidean projection onto $\mathcal{Q}^\circ$. Given a sequence of positive scalars $\sigma_k \uparrow \sigma_\infty \leqslant +\infty$ and a starting point $y^0 \in \mathcal{Q}^\circ$, the $(k+1)$-th iteration of the ALM is given by

$$\begin{cases} x^{k+1} \approx \arg\min \{f_k(x) := L_{\sigma_k}(x, y^k)\}, \\ y^{k+1} = \Pi_{\mathcal{Q}^\circ}[y^k + \sigma_k(\mathcal{B}x^{k+1} - b)], \quad k \geqslant 0. \end{cases} \tag{1}$$

The most appealing feature of the ALM for solving challenging optimization problems is its fast convergence rate. For the ALM applied to nonlinear programming (NLP), the classical results state that the generated dual sequence converges Q-linearly and the corresponding primal sequence converges R-linearly under the second order sufficient condition (SOSC), the linear independence constraint qualification and strict complementarity; see, e.g., [7, Propositions 2.7 & 3.2]. These assumptions automatically require the uniqueness of both the local optimal solution to the NLP problem and the corresponding multiplier. Various attempts have been made to relax these strong assumptions, such as in [31,12,13,22]. Among these works, Fernández and Solodov [22] showed that with a properly chosen initial multiplier and a sufficiently large penalty parameter in the case of an NLP problem, both the primal and dual sequences converge Q-linearly locally under solely the SOSC assumption. This nice result is made possible by the fact that the Karush-Kuhn-Tucker (KKT) solution mapping associated with the NLP problem is upper Lipschitz continuous if the SOSC is satisfied [18,34,32,38]. For convex NLP problems, the convergence rate of the ALM can also be derived through its connection with the dual proximal point algorithm (PPA) as championed by Rockafellar in [50]. Along this line, one can obtain the Q-linear convergence rate of the dual sequence generated by the ALM under the Lipschitz continuity of the dual solution mapping at the origin and certain stopping criteria on the inexact computation of the augmented Lagrangian subproblems [50, Proposition 3 & Theorem 2].

Despite the successes in the case of NLP, the aforementioned convergence rate results cannot be used to justify the necessity of employing the ALM in CCCP solvers. There are three main compelling reasons for this. Firstly, unlike the case of NLP, the Lipschitzian-like properties of the KKT solution mappings are difficult to be satisfied for CCCP when the cone $\mathcal{Q}$ is not a polyhedral set or the convex function $p$ is not piecewise linear-quadratic. For example, Bonnans and Shapiro constructed a convex quadratic SDP problem with a strongly convex objective function and a unique multiplier failing to possess an upper Lipschitz continuous KKT solution mapping [9, Example 4.54] (see also Example 1 in section 2.2 for a similar situation). Thus, directly applying the results in [22] to obtain the primal-dual convergence rate of the ALM for solving CCCP is not possible. Secondly, within the spirit of Rockafellar's work in [50], only the asymptotic Q-superlinear convergence rate of the dual sequence can be derived under the upper Lipschitz continuity of the dual solution mapping at the origin [37]. However, the lack of KKT residual information in this derived rate result poses a serious practical issue since the solution qualities in reliable solvers are almost always measured by the KKT residues. Thirdly, the subproblems in the ALM are often solved approximately in a CCCP solver. But the stopping criteria in [50], which involve the unknown optimal values of the subproblems, are difficult to be executed unless the objective function of the original problem is strongly convex. One plausible remedy for the latter two deficiencies encountered in the ALM is to instead adopt the proximal ALM, or in Rockafellar's terminology, the proximal method of multipliers [50]. However, a moment's thought would reveal that this is impractical as the primal-dual convergence rate of the proximal ALM would require a Lipschitzian-like property of the KKT solution mapping at the origin, which is too restrictive for CCCP as we have already mentioned above.

Now we are facing a dilemma as on the one hand, if the ALM is applied to solve CCCP, then only the Q-linear convergence of the dual sequence can be derived under the upper Lipschitz continuity of the dual solution mapping at the origin while on the other hand, if the proximal ALM is adopted, the assumption to ensure the linear convergence rate is too restrictive to hold even for the case when the original problem is strongly convex with a unique multiplier. This leads us to ask the following important question: *Is it possible for the KKT residues of the iterates generated by the ALM for solving CCCP to converge linearly without the Lipschitzian-like property of the KKT solution mapping at the origin?*

Before answering this question, we shall first conduct numerical experiments on the following convex least square SDP problem, which originates from [61] and [17, Example 3].

*Example 1* Consider the following SDP problem and its dual form:

$$\begin{array}{ll} \min_x \ \frac{1}{2}\|\mathcal{A}x - b\|^2 + \langle x, I \rangle & \max_{y,t} \ -\frac{1}{2}\|y\|^2 - \langle b, y \rangle - t \\ \text{s.t.} \ \langle E, x \rangle \leq 1, \ x \in \mathbb{S}_+^2; & \text{s.t.} \ \mathcal{A}^* y + tE + I \in \mathbb{S}_+^2, \ t \geq 0, \end{array}$$

where $\mathcal{A}x = B^{1/2}(x_{11}, x_{22})^T$ for all $x \in \mathbb{S}^2$ with $B = \begin{pmatrix} 3/2 & -2 \\ -2 & 3 \end{pmatrix}$, $b = B^{-1/2}(5/2, -1)^T$, $I$ is the $2 \times 2$ identity matrix and $E$ is the $2 \times 2$ matrix of all ones.

In Example 1, both the primal and the dual optimal solutions are unique, and the dual SOSC holds, but the KKT solution mapping is not Lipschitz continuous at the origin [17]. We apply the ALM to the primal form, with the subproblems being solved by the semismooth Newton method to the accuracy of $10^{-15}$. The following figure shows the KKT residual norm of the iterates, which is the maximum of the primal feasibility, dual feasibility and duality gap, against different choices of the penalty parameter $\sigma_k$. One can see that the KKT residues converge locally linearly; the linear convergence rate is a constant if $\sigma_k$ is fixed, and is decreasing if $\sigma_k$ is increased.

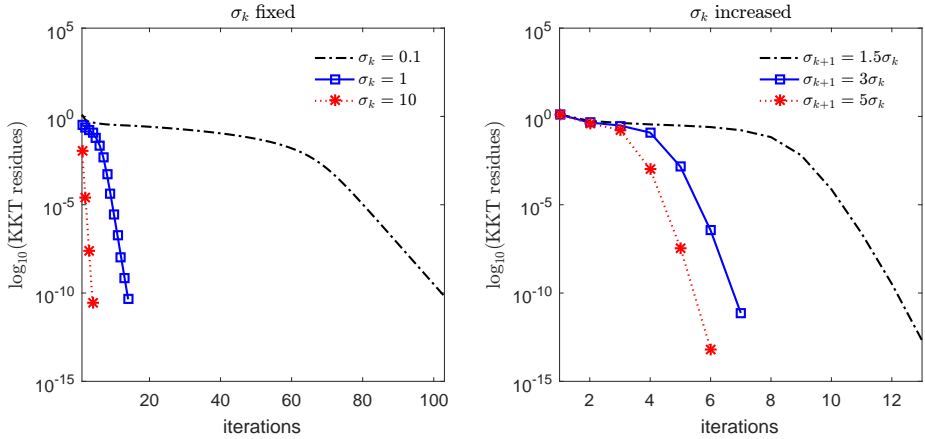

**Fig. 1** The KKT residual norm of the sequence generated by the ALM for solving Example 1 with different values of the penalty parameter $\sigma_k$.

The above numerical results shed some light on the possibility that the KKT residues of the iterates generated by the ALM could still converge linearly without the Lipschitz continuity of the KKT solution mapping at the origin. To establish a theorem of this property constitutes a major part of the present paper. Consequently, we fill the gap between the theoretical and practical performance of the ALM that has been missing so far. More specifically, our main contributions can be summarized below.

1. Under the stopping criteria proposed by Rockafellar in [50], we establish the R-linear convergence rate of the primal feasibility, complementarity and primal objective value under the quadratic growth condition on the dual problem. This quadratic growth condition can be satisfied even if the targeted problem admits multiple solutions and multipliers. The convergence rate becomes asymptotically superlinear if the penalty parameter tends to infinity.
2. Under the Robinson constraint qualification (RCQ), we propose fairly easy-to-implement stopping criteria for the inexact computation to the subproblems of the ALM to ensure the global convergence. More importantly, we prove that the KKT residues of the iterates generated by the ALM converge R-linearly under the quadratic growth condition on the dual problem. The convergence rate again becomes asymptotically superlinear if the penalty parameter tends to infinity.

We should mention that it is not completely new to guarantee the global convergence of the ALM for solving convex programming problems under implementable stopping criteria for the ALM subproblems. For example, Eckstein and Silva [21] proposed implementable relative stopping criteria for the augmented Lagrangian subproblems to ensure the global convergence of the dual iterates for convex

NLP; see also [2]. However, for the first time we prove, even with the subproblems being solved exactly, that the ALM can enjoy the much desired R-linear (asymptotically R-superlinear) convergence rate without the Lipschitzian-like property of the KKT solution mapping.

The remaining parts of this paper are organized as follows. In the next section, we provide some preliminary results on the convergence rates of the ALM. In section 3, we address the issue that the Lipschitzian-like property of the KKT solution mapping, with the presence of non-polyhedral set-valued mappings, is a much more restrictive assumption than that for either the primal or the dual solution mapping or even both. Section 4 is devoted to the study on the asymptotic R-superlinear convergence of the KKT residues for the iterates generated by the ALM for solving CCCP. Easy-to-implement stopping criteria are also provided in this section under the RCQ. The usefulness of the derived theoretical results are demonstrated for solving the convex quadratic SDP problems in section 5. Illustrative numerical experiments with real financial data are also conducted in this section. We end this paper with some concluding discussions in section 6.

**Notation.** We use $\mathbb{U}$, $\mathbb{V}$, $\mathbb{W}$, $\mathbb{X}$, $\mathbb{Y}$ and $\mathbb{Z}$ to denote finite dimensional real Euclidean spaces, $\mathbb{S}^n$ to denote the space of all $n \times n$ symmetric matrices, and $\mathbb{S}^n_+$ to denote the cone of all $n \times n$ symmetric positive semidefinite matrices. For any convex function $p : \mathbb{X} \to (-\infty, +\infty]$, we denote its effective domain as $\text{dom}\, p := \{x \in \mathbb{X} \mid p(x) < +\infty\}$, its epigraph as $\text{epi}\, p := \{(x,t) \in \mathbb{X} \times \mathbb{R} \mid p(x) \leq t\}$, its conjugate as $p^*(u) := \sup_{x \in \mathbb{X}}\{\langle x, u \rangle - p(x)\}$, $u \in \mathbb{X}$, and its proximal mapping as $\text{Prox}_p(x) := \arg\min_{u \in \mathbb{X}} \{\frac{1}{2}\|u - x\|^2 + p(u)\}$, $x \in \mathbb{X}$. We will often use the Moreau identity $x = \text{Prox}_p(x) + \text{Prox}_{p^*}(x)$ for any $x \in \mathbb{X}$ (c.f. [47, Theorem 31.5]). Let $D \subseteq \mathbb{X}$ be a convex set. Denote the indicator function over $D$ by $\delta_D(\cdot)$, i.e., for any $x \in \mathbb{X}$, $\delta_D(x) = 0$ if $x \in D$, and $\delta_D(x) = +\infty$ if $x \notin D$. We use $\text{cl}\,(D)$ and $\text{ri}\,(D)$ to denote the closure and relative interior of $D$, respectively. We write the distance of $x \in \mathbb{X}$ to $D$ by $\text{dist}(x, D) := \inf_{d \in D} \|d - x\|$. For a given closed convex set $D \subseteq \mathbb{X}$ and a given point $x \in \mathbb{X}$, the Euclidean projection of $x$ onto $D$ is denoted by $\Pi_D(x) := \arg\min\{\|x - d\| \mid d \in D\}$. For any $x \in D$, we use $\mathcal{T}_D(x)$ to denote the tangent cone of $D$ at $x$ and $\mathcal{N}_D(x)$ to denote the normal cone of $D$ at $x$. If $D$ is a closed convex cone, we use $D^\circ$ to denote the polar of $D$, i.e., $D^\circ := \{x \in \mathbb{X} \mid \langle x, d \rangle \leq 0, \ \forall d \in D\}$. For any set-valued mapping $\Gamma : \mathbb{U} \rightrightarrows \mathbb{V}$, we use $\text{gph}\,\Gamma$ to denote the graph of $\Gamma$, i.e., $\text{gph}\,\Gamma := \{(u, v) \in \mathbb{U} \times \mathbb{V} \mid v \in \Gamma(u)\}$. Let $\mathbb{B}_{\mathbb{U}}$ be the unit ball in $\mathbb{U}$ centered at the origin. For any $\bar{u} \in \mathbb{U}$ and $\varepsilon > 0$, denote $\mathbb{B}_\varepsilon(\bar{u}) := \{u \in \mathbb{U} \mid \|u - \bar{u}\| \leq \varepsilon\}$. We write the domain of $\Gamma$ as $\text{Dom}\,\Gamma := \{u \in \mathbb{U} \mid \Gamma(u) \neq \varnothing\}$.

## 2 Preliminary results on the convergence rates of the ALM

As mentioned in the introduction, the ALM can be taken as a dual application of the PPA for solving convex optimization problems [50]. Thus, one can obtain the convergence properties of the dual sequence generated by the ALM through the known results of the PPA. As a preparation for our subsequent study, we will review some known results and make necessary extensions along this line in this section.

### 2.1 The convergence rates of the PPA and ALM

Let $T : \mathbb{Z} \rightrightarrows \mathbb{Z}$ be a maximal monotone operator. Our aim is to find $z \in \mathbb{Z}$ such that

$$0 \in T(z). \tag{2}$$

Given a sequence of positive scalars $\sigma_k \uparrow \sigma_\infty \leq \infty$ and a starting point $z^0 \in \mathbb{Z}$, the $(k+1)$-th iteration of the PPA takes the form of

$$z^{k+1} \approx P_k(z^k) := (I + \sigma_k T)^{-1}(z^k), \quad \forall\, k \geq 0, \tag{3}$$

where $I$ is the identity operator in $\mathbb{Z}$. In [51], Rockafellar suggested the following criteria for computing $z^{k+1}$ approximately to ensure the global convergence and the convergence rate of the PPA:

(A) $\quad \|P_k(z^k) - z^{k+1}\| \leq \varepsilon_k, \quad \varepsilon_k \geq 0, \quad \sum_{k=1}^{\infty} \varepsilon_k < \infty,$

(B) $\quad \|P_k(z^k) - z^{k+1}\| \leq \eta_k \|z^{k+1} - z^k\|, \quad 0 \leq \eta_k < 1, \quad \sum_{k=1}^{\infty} \eta_k < \infty.$

In Rockafellar's original work [51], the asymptotic Q-superlinear convergence rate of $\{z^k\}$ is established under the Lipschitz continuity of $T^{-1}$ at the origin. Recall that a set-valued mappig $\Gamma : \mathbb{U} \to \mathbb{V}$ is called Lipschitz continuous at $\bar{u} \in \mathbb{U}$ if $\Gamma(\bar{u})$ admits a unique solution $\bar{v}$ and there exist positive constants $\kappa$ and $\varepsilon$ such that

$$\|v - \bar{v}\| \leq \kappa \|u - \bar{u}\|, \quad \forall\, v \in \Gamma(u), \quad \forall\, u \in \mathbb{B}_\varepsilon(\bar{u}).^1$$

Since the uniqueness assumption on the solution set to problem (2) may exclude many interesting instances, Luque in [37] made an important extension to Rockafellar's work by studying the linear convergence of the PPA under the following much more relaxed condition: there exist positive constants $\varepsilon$ and $\kappa$ such that

$$\operatorname{dist}(z, T^{-1}(0)) \leq \kappa \|u\|, \quad \forall\, z \in T^{-1}(u), \quad \forall\, u \in \mathbb{B}_\varepsilon(0).$$

In fact, this condition is exactly the so-called local upper Lipschitz continuity of $T^{-1}$ at the origin that was first coined by Robinson in [44]. A fundamental stability result of Robinson [45] states that every piecewise polyhedral mapping is locally upper Lipschitz continuous. Moreover, one key result in Sun's PhD thesis [57] says that a closed proper convex function is piecewise linear-quadratic if and only if its subdifferential is piecewise polyhedral (see also [52, Propositions 12.30 & 11.14]). Thus, when $T^{-1}$ is not piecewise polyhedral, such as $T^{-1}$ being the solution mapping of CCCP with $\mathcal{Q}$ being non-polyhedral or $p$ not being piecewise linear-quadratic, a weaker assumption on $T^{-1}$ may be needed for studying the convergence rate of the PPA. This leads us to focus on the calmness property. Recall that a set-valued mapping $\Gamma : \mathbb{U} \rightrightarrows \mathbb{V}$ is said to be calm at $\bar{u} \in \mathbb{U}$ for $\bar{v} \in \mathbb{V}$ (with modulus $\kappa$) (c.f. [19, 3.8(3H)]) if $(\bar{u}, \bar{v}) \in \operatorname{gph}\Gamma$ and there exist positive constants $\varepsilon$ and $\delta$ such that

$$\Gamma(u) \cap \mathbb{B}_\delta(\bar{v}) \subseteq \Gamma(\bar{u}) + \kappa \|u - \bar{u}\| \mathbb{B}_\mathbb{V}, \quad \forall\, u \in \mathbb{B}_\varepsilon(\bar{u}).$$

In the following proposition, we summarize and extend some useful results for the PPA developed in [51] and [37]. Among them, part (a) is from [51, Proposition 1(a)] and part (b) comes from [51, (2.11)]. The convergence rate result in part (c), whose proof is given in the appendix, is an extension of [51, Theorem 2] and [37, Theorem 2.1]. Note that for the case where the subproblems of the PPA are solved exactly, this relaxation has also been discussed by Leventhal in [35].

**Proposition 1** *Assume that $T^{-1}(0)$ is nonempty. Let $\{z^k\}$ be an infinite sequence generated by the PPA in (3) under criterion (A). Then the following three statements hold.*

*(a) $[z - P_k(z)]/\sigma_k \in T(P_k(z))$ for any $z \in \mathbb{Z}$ and $k \geq 0$.*

---

[1] The Lipschitz continuity of a set-valued mapping may refer to other properties elsewhere, such as in [52, Definition 9.26].

*(b) For any $\bar{z} \in T^{-1}(0)$, it holds that*

$$\|P_k(z^k) - \bar{z}\|^2 \leq \|z^k - \bar{z}\|^2 - \|P_k(z^k) - z^k\|^2, \quad \forall\, k \geq 0.$$

*(c) The whole sequence $\{z^k\}$ converges to some $z^\infty \in T^{-1}(0)$. If in the PPA, the stopping criterion (B) is also employed and the mapping $T^{-1}$ is calm at the origin for $z^\infty$ with modulus $\kappa$, then there exists $\bar{k} \geq 0$ such that for all $k \geq \bar{k}$,*

$$\mathrm{dist}\,(z^{k+1}, T^{-1}(0)) \leq \mu_k \,\mathrm{dist}\,(z^k, T^{-1}(0)),$$

*where $\mu_k := \left[\eta_k + (\eta_k + 1)\kappa/\sqrt{\kappa^2 + \sigma_k^2}\,\right]/(1 - \eta_k) \to \mu_\infty := \kappa/\sqrt{\kappa^2 + \sigma_\infty^2}$ ($\mu_\infty = 0$ if $\sigma_\infty = +\infty$).*

Let us now move on to the ALM. To proceed, we first introduce some notation. Let $l : \mathbb{X} \times \mathbb{Y} \to [-\infty, +\infty]$ be the Lagrangian function of (P) in the extended form:

$$l(x, y) := \begin{cases} f^0(x) + \langle y, \mathcal{B}x - b\rangle & x \in \mathrm{dom}\, f^0,\ y \in \mathcal{Q}^\circ, \\ -\infty & x \in \mathrm{dom}\, f^0,\ y \notin \mathcal{Q}^\circ, \\ +\infty & x \notin \mathrm{dom}\, f^0. \end{cases}$$

The Lagrangian dual of (P) takes the form of

$$\max_{y \in \mathbb{Y}}\ \{g^0(y) := \inf_{x \in \mathbb{X}} l(x, y)\}, \quad \text{s.t. } y \in \mathcal{Q}^\circ. \tag{D}$$

The essential objective functions of (P) and (D) are given by

$$f(x) := \sup_{y \in \mathbb{Y}} l(x, y) = \begin{cases} f^0(x) & \mathcal{B}x \in b + \mathcal{Q}, \\ +\infty & \text{otherwise}; \end{cases} \qquad g(y) := \inf_{x \in \mathbb{X}} l(x, y) = \begin{cases} g^0(y) & y \in \mathcal{Q}^\circ, \\ -\infty & \text{otherwise}. \end{cases}$$

Denote the mappings $T_l : \mathbb{X} \times \mathbb{Y} \rightrightarrows \mathbb{X} \times \mathbb{Y}$, $T_f : \mathbb{X} \rightrightarrows \mathbb{X}$ and $T_g : \mathbb{Y} \rightrightarrows \mathbb{Y}$ by

$$T_l(x, y) := \{(u, v) \in \mathbb{X} \times \mathbb{Y} \mid (u, -v) \in \partial l(x, y)\},\ (x, y) \in \mathbb{X} \times \mathbb{Y},\ T_f := \partial f,\ T_g := -\partial g, \tag{4}$$

where $\partial l$ is the subdifferential of the convex-concave function $l$ and $\partial g$ is the subdifferential of the concave function $g$; see [47, Sections 30 & 35] for the definitions of such subdifferentials. All of the three mappings $T_l$, $T_f$ and $T_g$ are maximal monotone operators under the settings of (P) [47, Corollaries 37.5.2 & 31.5.2]. As explained in [50], the inverse of these mappings can be taken as the solution mappings of the corresponding perturbed problems. Specifically, consider the following linearly perturbed form of problem (P):

$$\begin{aligned} \min\ & f^0(x) - \langle x, u\rangle \\ \text{s.t.}\ & \mathcal{B}x + v \in b + \mathcal{Q}, \end{aligned} \tag{P($u,v$)}$$

where $(u, v) \in \mathbb{X} \times \mathbb{Y}$ are perturbation parameters. Then

$$\begin{cases} T_l^{-1}(u, v) = \text{set of all KKT points to (P($u,v$))}, \\ T_f^{-1}(u) = \text{set of all optimal solutions to (P($u,0$))}, \\ T_g^{-1}(v) = \text{set of all optimal solutions to (D($0,v$))}, \end{cases}$$

where $(\mathrm{D}(u,v))$ is the ordinary dual of $(\mathrm{P}(u,v))$ for any $(u,v) \in \mathbb{X} \times \mathbb{Y}$. Thus, we call $T_l^{-1}$ the KKT solution mapping, $T_f^{-1}$ the primal solution mapping, and $T_g^{-1}$ the dual solution mapping, all with respect to (P).

For the ALM, the following criteria on the approximate computation of $x^{k+1}$ in (1) are considered by Rockafellar in [50]:

(A′) $\quad f_k(x^{k+1}) - \inf f_k \leqslant \varepsilon_k^2/2\sigma_k, \quad \varepsilon_k \geqslant 0, \quad \sum_{k=0}^{\infty} \varepsilon_k < \infty,$

(B′) $\quad f_k(x^{k+1}) - \inf f_k \leqslant (\eta_k^2/2\sigma_k)\|y^{k+1} - y^k\|^2, \quad 0 \leqslant \eta_k < 1, \quad \sum_{k=0}^{\infty} \eta_k < \infty,$

($\widetilde{B}$) $\quad \mathrm{dist}\,(0, \partial f_k(x^{k+1})) \leqslant (\eta'_k/\sigma_k)\|y^{k+1} - y^k\|, \quad 0 \leqslant \eta'_k \to 0.$

The connection between the ALM applied to (P) and the PPA applied to (D) is revealed in the next lemma, where part (a) is taken from [50, Proposition 6] and part (b) is taken from [50, (4.21)].

**Lemma 1** *For any $k \geqslant 0$ and $x^{k+1} \in \mathbb{X}$, let $P_k = (I + \sigma_k T_g)^{-1}$ and $y^{k+1} = \Pi_{Q^\circ}[y^k + \sigma_k(\mathcal{B}x^{k+1} - b)]$. Then the following inequalities hold:*

*(a) $(1/2\sigma_k)\|y^{k+1} - P_k(y_k)\|^2 \leqslant f_k(x^{k+1}) - \inf f_k$.*

*(b) $\mathrm{dist}\,(0, T_l^{-1}(x^{k+1}, y^{k+1})) \leqslant [\mathrm{dist}^2(0, \partial f_k(x^{k+1})) + (1/\sigma_k^2)\|y^{k+1} - y^k\|^2]^{1/2}$.*

It can be seen from Lemma 1 that, when using the ALM for (P) under criteria (A′) and (B′), we are, in effect, executing the PPA for $T_g = -\partial g$ under criteria (A) and (B). This fact leads to the following two propositions regarding the global convergence and local convergence rate of the ALM.

**Proposition 2** *Assume that $T_g^{-1}(0)$ is nonempty. Let $\{(x^k, y^k)\}$ be an infinite sequence generated by the ALM for (P) under criterion (A′). Then, the whole sequence $\{y^k\}$ is bounded and converges to some $y^\infty \in T_g^{-1}(0)$, and the sequence $\{x^k\}$ satisfies for all $k \geqslant 0$,*

$$\|\mathcal{B}x^{k+1} - b - \Pi_{\mathcal{Q}}[\mathcal{B}x^{k+1} - b + y^k/\sigma^k]\| = (1/\sigma_k)\|y^{k+1} - y^k\| \to 0, \tag{5a}$$

$$f^0(x^{k+1}) - \inf(\mathrm{P}) \leqslant f_k(x^{k+1}) - \inf f_k + (1/2\sigma_k)(\|y^k\|^2 - \|y^{k+1}\|^2). \tag{5b}$$

*Moreover, if (P) admits a nonempty and bounded solution set, then the sequence $\{x^k\}$ is also bounded, and all of its accumulation points are optimal solutions to (P).*

The above proposition is essentially adopted from [50, Theorem 4]. Note that the inequalities (5a) and (5b) are slightly different from [50, Theorem 4 (4.13) & (4.14)]. We directly take the inequalities (4.4) and (4.18) in the proof of [50, Theorem 4] to serve the purpose of our later developments.

**Proposition 3** *Let $\{(x^k, y^k)\}$ be an infinite sequence generated by the ALM for (P) under criterion (A′), for which $\{y^k\}$ converges to $y^\infty$.*

*(a) If $T_g^{-1}$ is calm at the origin for $y^\infty$ with modulus $\kappa_g$, then under criterion (B′), there exists $\bar{k} \geqslant 0$ such that for all $k \geqslant \bar{k}$,*

$$\mathrm{dist}\,(y^{k+1}, T_g^{-1}(0)) \leqslant \mu_k \mathrm{dist}\,(y^k, T_g^{-1}(0)), \tag{6}$$

*where $\mu_k := \left[\eta_k + (\eta_k + 1)\kappa_g/\sqrt{\kappa_g^2 + \sigma_k^2}\right]/(1 - \eta_k) \to \mu_\infty := \kappa_g/\sqrt{\kappa_g^2 + \sigma_\infty^2}$ ($\mu_\infty = 0$ if $\sigma_\infty = +\infty$).*

*(b) If in addition to $(B')$ and the calmness condition on $T_g^{-1}$, one has criterion $(\widetilde{B})$ and $T_l^{-1}$ is upper Lipschitz continuous at the origin with modulus $\kappa_l$, then there exists $\tilde{k} \geq 0$ such that for all $k \geq \tilde{k}$,*

$$\mathrm{dist}\,(x^{k+1}, T_f^{-1}(0)) \leq \mu_k' \|y^{k+1} - y^k\|,$$

*where $\mu_k' := (\kappa_l/\sigma_k)(1 + \eta_k'^2) \to \mu_\infty' := \kappa_l/\sigma_\infty$ ($\mu_\infty' = 0$ if $\sigma_\infty = +\infty$).*

Part (a) of the above proposition is a consequence of criteria $(A')$ and $(B')$, Proposition 1(c) and Lemma 1(a); part (b), whose proof is given in the appendix, is an extension of [50, Theorem 5] by relaxing the Lipschitz continuity of $T_l^{-1}$ at the origin, which is too restrictive in our context. Proposition 3 provides the asymptotic Q-superlinear convergence rate of $\{y^k\}$ under the calmness of the dual solution mapping $T_g^{-1}$ at the origin. In fact, by applying Lemma 3 that will be given in section 4, one can also derive the asymptotic R-superlinear rate of $\{x^k\}$ under the assumptions in part (b). However, the assumption on the upper Lipschitz continuity of $T_l^{-1}$ is likely to fail for CCCP considered in this paper. This phenomenon can be seen clearly later in section 3.

2.2 The quadratic growth condition of CCCP

Proposition 3 says that the asymptotic Q-superlinear convergence of the dual sequence generated by the ALM holds under the calmness assumption of $T_g^{-1}$ at the origin. Then one may wonder to what extent this condition is satisfied for CCCP, especially when $T_g$ is a non-polyhedral set-valued mapping. In this subsection, we attempt to characterize the calmness of the dual solution mapping by relating it to the known results in the existing literature.

Since the function $h$ is assumed to be essentially smooth, it is known from [47, Theorem 26.1] that $\partial h(x)$ is nonempty if and only if $x \in \mathrm{int}\,(\mathrm{dom}\,h)$, where in fact $\partial h(x)$ consists of $\nabla h(x)$ alone. Thus, the KKT optimality condition of (P) can be written the form of

$$0 \in \mathcal{A}^* \nabla h(\mathcal{A}x) + c + \partial p(x) + \mathcal{B}^* y, \quad y \in \mathcal{N}_\mathcal{Q}(\mathcal{B}x - b), \quad (x,y) \in \mathbb{X} \times \mathbb{Y}. \tag{7}$$

Throughout this paper, we always assume that the KKT system (7) has at least one solution. Denote $\mathrm{SOL}_\mathrm{P}$ and $\mathrm{SOL}_\mathrm{D}$ as the sets of all the optimal solutions to (P) and (D), respectively. Then $\bar{x} \in \mathrm{SOL}_\mathrm{P}$ and $\bar{y} \in \mathrm{SOL}_\mathrm{D}$ if and only if $(\bar{x}, \bar{y})$ solves (7) [47, Theorem 28.3]. Denote $\mathrm{F}_\mathrm{P}$ and $\mathrm{F}_\mathrm{D}$ as the set of all the primal and dual feasible solutions, i.e.,

$$\mathrm{F}_\mathrm{P} := \{x \in \mathbb{X} \mid \mathcal{B}x \in b + \mathcal{Q},\ x \in \mathrm{dom}\,f^0\}, \quad \mathrm{F}_\mathrm{D} := \{y \in \mathbb{Y} \mid y \in \mathrm{dom}\,(-g^0) \cap \mathcal{Q}^\circ\}.$$

The quadratic growth condition for (P) at $\bar{x} \in \mathrm{SOL}_\mathrm{P}$ is said to hold if there exist positive constants $\kappa_1 > 0$ and $\varepsilon_1 > 0$ such that

$$f^0(x) \geq \inf(\mathrm{P}) + \kappa_1 \mathrm{dist}^2(x, \mathrm{SOL}_\mathrm{P}), \quad \forall\, x \in \mathrm{F}_\mathrm{P} \cap \mathbb{B}_{\varepsilon_1}(\bar{x}). \tag{8}$$

The quadratic growth condition for (D) at $\bar{y} \in \mathrm{SOL}_\mathrm{D}$ is said to hold if there exist positive constants $\kappa_2 > 0$ and $\varepsilon_2 > 0$ such that

$$-g^0(y) \geq -\sup(\mathrm{D}) + \kappa_2 \mathrm{dist}^2(y, \mathrm{SOL}_\mathrm{D}), \quad \forall\, y \in \mathrm{F}_\mathrm{D} \cap \mathbb{B}_{\varepsilon_2}(\bar{y}). \tag{9}$$

The constants $\kappa_1$ and $\kappa_2$ are called the quadratic growth modulus for (P) at $\bar{x}$ and for (D) at $\bar{y}$, respectively. It can be derived from [3, Theorem 3.3] that the calmness of the dual solution mapping at the origin for a dual optimal solution is equivalent to the quadratic growth condition at that optimal solution. This is made precise in the proposition below.

**Proposition 4** *Assume that* $\mathrm{SOL}_D$ *is nonempty. Let* $\bar{y} \in \mathrm{SOL}_D$. *The following statements are equivalent to each other:*
*(i) The mapping $T_g^{-1}$ is calm at the origin for $\bar{y}$.*
*(ii) The quadratic growth condition (9) for (D) holds at $\bar{y}$.*

*More specifically, if (9) holds with quadratic growth modulus $\kappa$, then $T_g^{-1}$ is calm at the origin for $\bar{y}$ with modulus $1/\kappa$; conversely, if $T_g^{-1}$ is calm at the origin for $\bar{y}$ with modulus $\kappa'$, then (9) holds for any $\kappa \in (0, 1/(4\kappa'))$.*

*Proof* By [3, Theorem 3.3] and the fact that $g(y) = g^0(y)$ for any $y \in \mathrm{F}_D$, the quadratic growth condition for (D) holds at $\bar{y} \in \mathrm{SOL}_D$ if and only if $-\partial g = \partial(-g)$ is metrically subregular (see [19, 3.8(3H)] for its definition) at $\bar{y}$ for the origin. The latter property is the same as the calmness of the mapping $T_g^{-1} = (-\partial g)^{-1}$ at the origin for $\bar{y}$ with the same modulus, as stated in [19, Theorem 3H.3]. □

The above proposition allows us to characterize the calmness of the mapping $T_g^{-1}$ of CCCP via the dual quadratic growth conditions studied in the recent work [14]. However, the lack of an explicit expression of the dual objective function makes it difficult to apply the known results directly. To see this, let the function $\phi : \mathbb{X} \to (-\infty, +\infty]$ be given by $\phi(x) := h(\mathcal{A}x) + p(x)$ for $x \in \mathbb{X}$. Direct computations show that
$$g^0(y) = -\phi^*(-\mathcal{B}^*y - c) - \langle b, y \rangle, \quad y \in \mathbb{Y},$$
where the conjugate function $\phi^*$ takes the form of
$$\phi^*(v) = \inf_{w \in \mathbb{W}} \{h^*(w) + p^*(v - \mathcal{A}^*w)\}, \quad v \in \mathbb{X}.$$

To alleviate the aforementioned difficulty, we introduce the following auxiliary problem for (D):
$$\max_{(w,y,s) \in \mathbb{W} \times \mathbb{Y} \times \mathbb{X}} g^0(w, y, s) := -h^*(w) - \langle b, y \rangle - p^*(s) \quad \text{(D2)}$$
$$\text{s.t.} \quad \mathcal{A}^*w + \mathcal{B}^*y + s + c = 0, \quad y \in \mathcal{Q}^\circ.$$

Note that we used the same notation $g^0$ to denote the dual objective function of (D2) but there is no danger of confusion since it involves three variables $(w, y, s)$. Similarly to the notation $\mathrm{SOL}_D$, we use $\mathrm{SOL}_{D2}$ to denote the set of all optimal solutions to (D2). Also denote by $\mathrm{F}_{D2}$ the set of all feasible solutions to (D2):
$$\mathrm{F}_{D2} := \{(w, y, s) \in \mathrm{dom}\, h^* \times \mathbb{Y} \times \mathbb{X} \mid \mathcal{A}^*w + \mathcal{B}^*y + s + c = 0, \, y \in \mathcal{Q}^\circ, \, s \in \mathrm{dom}\, p^*\}.$$

In fact, the KKT optimality condition in terms of (D2) can be written, for $(w, y, s, x) \in \mathbb{W} \times \mathbb{Y} \times \mathbb{X} \times \mathbb{X}$, as
$$\mathcal{A}x \in \partial h^*(w), \quad \mathcal{B}x - b \in \mathcal{N}_{\mathcal{Q}^\circ}(y), \quad x \in \partial p^*(s), \quad \mathcal{A}^*w + \mathcal{B}^*y + s + c = 0. \quad (10)$$

The quadratic growth condition for (D2) at $(\bar{w}, \bar{y}, \bar{s}) \in \mathrm{SOL}_{D2}$ is said to hold if there exist positive constants $\kappa_3 > 0$ and $\varepsilon_3 > 0$ such that
$$-g^0(w, y, s) \geq -\sup(\mathrm{D2}) + \kappa_3 \, \mathrm{dist}^2((w, y, s), \mathrm{SOL}_{D2}), \quad \forall \, (w, y, s) \in \mathrm{F}_{D2} \cap \mathbb{B}_{\varepsilon_3}(\bar{w}, \bar{y}, \bar{s}). \quad (11)$$

It is known from [47, Theorem 26.3] that if $h$ is essentially smooth, then its conjugate function $h^*$ is essentially strictly convex, implying that the vector $w$ is unique over $\mathrm{SOL}_{D2}$. For notational simplicity, we write this vector as $\bar{w}$. In the following lemma, we establish the equivalence between the quadratic growth conditions for problems (D) and (D2), whose proof can be found in the appendix.

**Lemma 2** *Assume that $h^*$ is locally strongly convex on $\mathrm{dom}\,h^*$. Then the quadratic growth condition for (D) holds at $\bar{y} \in \mathrm{SOL}_D$ if and only if the quadratic growth condition for (D2) holds at $(\bar{w}, \bar{y}, -\mathcal{A}^*\bar{w} - \mathcal{B}^*\bar{y} - c) \in \mathrm{SOL}_{D2}$.*

Equipped with the preparations in Proposition 4 and Lemma 2, one can obtain sufficient conditions for the calmness of the dual solution mapping associated with a rich class of non-polyhedral CCCP via the recent results established in [14]. In order not to deviate too much from the main purpose of this paper, we would not repeat those results here. Instead, we take convex quadratic SDP problems as an example to illustrate these sufficient conditions in section 5.

## 3 The Lipschitzian-like property of the KKT solution mapping can easily fail

As mentioned in the previous section, from Proposition 3 and Lemma 3 (in section 4), one can derive that the primal sequence generated by the ALM converges asymptotically R-superlinearly if the KKT solution mapping $T_l^{-1}$ is assumed to be upper Lipschitz continuous at the origin. However, this is a restrictive assumption for non-polyhedral CCCP. As what have already been indicated in Example 1, even if the dual SOSC holds and a unique KKT point is admitted, the mapping $T_l^{-1}$ can still fail to be upper Lipschitz continuous at the origin. This is completely different from the case of NLP, for which the SOSC implies the upper Lipschitz continuity of the KKT solution mapping at the origin, see, e.g., [18, 34, 32, 38]. In this section, we shall show that even a weaker assumption on CCCP – the calmness of the KKT solution mapping, can still be difficult to be satisfied.

In what follows, we shall first consider the case that the KKT system (7) admits a unique solution, in which case the calmness property coincides with the isolated calmness property. Recall that a set-valued mapping $\Gamma : \mathbb{U} \rightrightarrows \mathbb{V}$ is said to be isolated calm at $\bar{u} \in \mathbb{U}$ for $\bar{v} \in \mathbb{V}$ (c.f. [19, Section 3.9(3I)]) if $(\bar{u}, \bar{v}) \in \mathrm{gph}\,\Gamma$ and there exist positive constants $\kappa$, $\varepsilon$ and $\delta$ such that

$$\|v - \bar{v}\| \leqslant \kappa \|u - \bar{u}\|, \quad \forall\, v \in \Gamma(u) \cap \mathbb{B}_\delta(\bar{v}), \quad u \in \mathbb{B}_\varepsilon(\bar{u}). \tag{12}$$

If in addition, the mapping $\Gamma$ is also locally nonempty-valued, i.e., $\Gamma(u) \cap \mathbb{B}_\delta(\bar{v}) \neq \varnothing$ for any $u \in \mathbb{B}_\varepsilon(\bar{u})$, the set-valued mapping $\Gamma$ is called robustly isolated calm at $\bar{u}$ for $\bar{v}$ [17].

Based on the recent work in [17, Theorem 24], one can obtain the following characterization on the robust isolated calmness of the KKT solution mapping at the origin when the closed convex cone $\mathcal{Q}$ is $\mathcal{C}^2$-cone reducible and the function $p$ is the indicator function over a $\mathcal{C}^2$-cone reducible set ([9, Definition 3.135]). It is worth mentioning that the class of $\mathcal{C}^2$-cone reducible sets is rich, and it includes polyhedral convex sets, second order cones, the cones of symmetric and positive semidefinite matrices, the epigraph of the nuclear norm, and their Cartesian products [54, 16].

**Proposition 5** *Assume that the cone $\mathcal{Q}$ is $\mathcal{C}^2$-cone reducible and the function $p$ in (P) is $\delta_\mathcal{K}(\cdot)$, the indicator function over a $\mathcal{C}^2$-cone reducible set $\mathcal{K} \subseteq \mathbb{X}$. Let $(\bar{x}, \bar{y})$ be a KKT solution of (P). Then the KKT solution mapping $T_l^{-1}$ is isolated calmness at the origin for $(\bar{x}, \bar{y})$ if and only if it is isolated calmness at the origin for $(\bar{x}, \bar{y})$, which is also equivalent to have both the SOSC for (P) at $\bar{x}$ and the following strict Robinson constraint qualification for (P) at $\bar{x}$ for $\bar{y}$:*

$$\begin{pmatrix} \mathcal{B} \\ \mathcal{I} \end{pmatrix} \mathbb{X} + \begin{pmatrix} \mathcal{T}_\mathcal{Q}(\mathcal{B}\bar{x} - b) \cap \bar{y}^\perp \\ \mathcal{T}_\mathcal{K}(\bar{x}) \cap (\mathcal{A}^*\nabla h(\mathcal{A}\bar{x}) + \mathcal{B}^*\bar{y} + c)^\perp \end{pmatrix} = \begin{pmatrix} \mathbb{Y} \\ \mathbb{X} \end{pmatrix}.$$

One may refer to [15] for an analogous result to Proposition 5 for the case when the function $p$ in (P) is the nuclear norm function defined over $\mathbb{X} = \mathbb{R}^{m \times n}$ and $\mathcal{Q}$ is a convex polyhedral cone.

In fact, we can also characterize the Lipschitz continuity of the KKT solution mapping based on Proposition 5. The key factor to this characterization is the equivalence between the Lipschitz continuity and the robust isolated calmness of any maximal monotone mapping, which is given by the proposition below.

**Proposition 6** *Let $\Gamma : \mathbb{U} \to \mathbb{U}$ be a maximal monotone mapping. For any $(\bar{u}, \bar{v}) \in \text{gph}\,\Gamma$, the mapping $\Gamma$ is Lipschitz continuous at $\bar{u}$ if and only if it is robustly isolated calm at $\bar{u}$ for $\bar{v}$.*

*Proof* First let us prove the "if" part. By the definitions of robust isolated calmness and Lipschitz continuity, it suffices to show that $\Gamma$ is locally nonempty-valued at $\bar{u}$ for $\bar{v}$. This can be obtained by a result of Rockafellar that $\bar{u} \in \text{int}(\text{Dom}\,\Gamma)$ if $\Gamma$ is locally uniformly bounded [48, Thoerem 1], and the latter property is guaranteed by the definition of the Lipschitz continuity of $\Gamma$ at $\bar{u}$.

Now we prove the "only if" part. The assumed isolated calmness of $\Gamma$ at $\bar{u}$ for $\bar{v}$ implies the existence of positive constants $\varepsilon$, $\delta$ and $\kappa$ such that

$$\|v - \bar{v}\| \leq \kappa \|u - \bar{u}\|, \quad \forall\, v \in \Gamma(u) \cap \mathbb{B}_\delta(\bar{v}), \quad u \in \mathbb{B}_\varepsilon(\bar{u}). \tag{13}$$

Without loss of generality, let us assume that $\varepsilon < \delta/\kappa$. Suppose on the contrary that $\Gamma$ is not Lipschitz continuous at $\bar{u}$. By shrinking $\varepsilon$ if necessary, we may choose a sequence $\{(u^j, v^j)\}_{j \geq 1}$ satisfying

$$u^j \in \mathbb{B}_\varepsilon(\bar{u}) \setminus \{\bar{u}\}, \quad v^j \in \Gamma(u^j), \quad \|v^j - \bar{v}\| \geq t_j \|u^j - \bar{u}\| \text{ for some } t_j \to \infty. \tag{14}$$

We may also assume that $t_j > \kappa$ for any $j \geq 1$. The inequalities in (13) and (14) together imply that $v^j \notin \mathbb{B}_\delta(\bar{v})$ for any $j \geq 1$, or equivalently, $\|v^j - \bar{v}\| > \delta$. Consider a fixed but arbitrary index $j \geq 1$. Since $\Gamma$ is assumed to be locally nonempty, we know that there exists $\tilde{v}^j \in \Gamma(u^j) \cap \mathbb{B}_\delta(\bar{v})$. Denote the constant $\gamma_j := \frac{1}{2}(\kappa \|u^j - \bar{u}\| + \delta)$ and the function $\xi_j : \mathbb{R} \to \mathbb{R}$ as $\xi_j(\lambda) = \|\lambda v^j + (1-\lambda)\tilde{v}^j - \bar{v}\|$ for $\lambda \in \mathbb{R}$. Note that $\gamma_j \in (\kappa \|u^j - \bar{u}\|, \delta)$ since $u^j \in \mathbb{B}_\varepsilon(\bar{u}) \subsetneq \mathbb{B}_{\delta/\kappa}(\bar{u})$. Obviously $\xi_j(\cdot)$ is continuous on $\mathbb{R}$ with $\xi_j(0) < \gamma_j$ and $\xi_j(1) > \gamma_j$. Then by the intermediate value theorem, we get the existence of $\lambda_j \in (0, 1)$ for which $\xi_j(\lambda_j) = \gamma_j$. Since $\Gamma$ is a maximal monotone mapping, it is easy to check that $\lambda_j v^j + (1-\lambda_j)\tilde{v}^j \in \Gamma(u^j)$. Then $\lambda_j v^j + (1-\lambda_j)\tilde{v}^j \in \Gamma(u^j) \cap \mathbb{B}_\delta(\bar{v})$ since $\xi_j(\lambda_j) = \gamma_j < \delta$. On the other hand,

$$\|\lambda_j v^j + (1-\lambda_j)\tilde{v}^j - \bar{v}\| = \xi_j(\lambda_j) = \gamma_j > \kappa \|u^j - \bar{u}\|.$$

Thus, we get a contradiction to (13), which implies that $\Gamma$ must be Lipschitz continuous at $\bar{u}$. □

Propositions 5 and 6 reveal the reason behind the failure of the Lipschitz continuity of $T_l^{-1}$ at the origin in Example 1, that is, the lack of the strict Robinson constraint qualification. Different from the case of NLP, there is a gap between such a constraint qualification (known as the strict Mangasarian-Fromovitz constraint qualification in the case of NLP) and the uniqueness of the multiplier.

Now we move on to the case that the KKT system (7) admits multiple solutions. As indicated in [14, Theorem 16 & Proposition 17], to guarantee the calmness of $T_g^{-1}$ at the origin for a dual optimal solution, it suffices for (P) to possess a partial strict complementarity solution with respect to the non-polyhedral complementarity system. However, to ensure the calmness of the mapping $T_l^{-1}$ at the origin for a KKT point is much harder, as can be seen from the following example.

*Example 2* Consider the following SDP problem and its dual:

$$\min_{(x_1,x_2)\in\mathbb{S}^2\times\mathbb{R}} \delta_{\mathbb{S}^2_+}(x_1) \qquad\qquad \max_{s\in\mathbb{S}^2} s_{22} - \delta_{\mathbb{S}^2_-}(s)$$
$$\text{s.t.} \quad x_1 + \begin{pmatrix} 0 & x_2 \\ x_2 & -x_2 \end{pmatrix} = \begin{pmatrix} 0 & 0 \\ 0 & 1 \end{pmatrix}; \qquad\qquad \text{s.t.} \ 2s_{12} - s_{22} = 0.$$

It is easy to check that

$$\text{SOL}_\text{P} = \left\{\bar{x}_1 = \begin{pmatrix} 0 & 0 \\ 0 & 1 \end{pmatrix}, \bar{x}_2 = 0\right\}, \quad \text{SOL}_\text{D} = \left\{\bar{s} = \begin{pmatrix} t & 0 \\ 0 & 0 \end{pmatrix} \ \bigg|\ t \leq 0\right\}.$$

For any $(\bar{x}_1, \bar{x}_2, \bar{s}) \in \text{SOL}_\text{P} \times \text{SOL}_\text{D}$ with $\bar{s}_{11} < 0$, it holds that $\text{rank}(\bar{x}_1) + \text{rank}(\bar{s}) = 2$, or equivalently, $\bar{s} \in \text{ri}(\partial \delta_{\mathbb{S}^2_+}(\bar{x}_1))$. Hence, the dual solution mapping $T_g^{-1}$ is calm at the origin for any $\hat{s} \in \text{SOL}_\text{D}$ ($\hat{s}$ may be different from $\bar{s}$) [14, Theorem 16 & Proposition 17]. One may also obtain the calmness of $T_f^{-1}$ at the origin for any $(\bar{x}_1, \bar{x}_2) \in \text{SOL}_\text{P}$ as the SOSC for the primal problem is satisfied [14, Theorem 12]. However, the KKT solution mapping $T_l^{-1}$ fails to be calm at the origin for $(\bar{x}_1, \bar{x}_2, \bar{s}') \in T_l^{-1}(0)$ with $\bar{s}'_{11} = 0$. This can be seen as follows: for any $\epsilon > 0$, consider the perturbation parameters $u(\epsilon) := \left(\begin{pmatrix} \epsilon & \epsilon \\ \epsilon & 2\epsilon \end{pmatrix}, 0\right) \in \mathbb{S}^2 \times \mathbb{R}$ and $v(\epsilon) := \begin{pmatrix} -\epsilon & 0 \\ 0 & 0 \end{pmatrix} \in \mathbb{S}^2$. Then one can show from the KKT optimality condition that

$$(x_1(\epsilon), x_2(\epsilon), s(\epsilon)) := \left(\begin{pmatrix} \epsilon & -\sqrt{\epsilon} \\ -\sqrt{\epsilon} & 1+\sqrt{\epsilon} \end{pmatrix}, \sqrt{\epsilon}, \begin{pmatrix} \epsilon & \epsilon \\ \epsilon & 2\epsilon \end{pmatrix}\right) \in T_l^{-1}(u(\epsilon), v(\epsilon)).$$

Also one can readily verify that

$$\|(u(\epsilon), v(\epsilon))\| = 2\sqrt{2}\epsilon, \quad \text{dist}(s(\epsilon), \text{SOL}_\text{D}) = \left\|\begin{pmatrix} 0 & \epsilon \\ \epsilon & 2\epsilon \end{pmatrix}\right\| = \sqrt{6}\epsilon,$$
$$\text{dist}((x_1(\epsilon), x_2(\epsilon)), \text{SOL}_\text{P}) = \left((\sqrt{\epsilon})^2 + \left\|\begin{pmatrix} \epsilon & -\sqrt{\epsilon} \\ -\sqrt{\epsilon} & \sqrt{\epsilon} \end{pmatrix}\right\|^2\right)^{1/2} \geq 2\sqrt{\epsilon}.$$

Thus, there cannot exist positive constants $\kappa$ and $\delta$ such that

$$\text{dist}((x_1, x_2, s), T_l^{-1}(0)) \leq \kappa \|(u, v)\|, \quad \forall\ (x_1, x_2, s) \in T_l^{-1}(u, v) \cap \mathbb{B}_\delta(\bar{x}_1, \bar{x}_2, \bar{s}'),$$

showing that $T_l^{-1}$ cannot be calm at the origin for $(\bar{x}_1, \bar{x}_2, \bar{s}') \in T_l^{-1}(0)$ with $\bar{s}'_{11} = 0$.

Example 2 says that the mapping $T_l^{-1}$ associated with the non-polyhedral CCCP may fail to possess the calmness property at the origin for a KKT point even if both $T_f^{-1}$ and $T_g^{-1}$ are calm at the origin for a corresponding optimal solution. Therefore, additional conditions must be imposed in order to guarantee the Lipschitz continuity and calmness of $T_l^{-1}$ at the origin for CCCP.

We shall end this section with the following remarks. Though the uniqueness of the solution and multiplier, strict complementarity, SOSC, and constraint nondegeneracy[2] are all generic properties for many structured conic programming problems, including the linear and nonlinear SDP [1,53,40,20], some of these properties often fail to hold for problems arising from various interesting applications. In general, the quadratic growth condition for (P)/(D) is much weaker than the Lipschitz continuity/calmness of the KKT solution mapping at the origin. While the former condition is satisfied under either the SOSC for (P)/(D) or the existence of a partial strict complementarity solution [14], the latter may fail even if both the quadratic growth conditions for (P) and (D) hold.

---

[2] For NLP, the constraint nondegeneracy coincides with the linear independence constraint qualification [46].

## 4 The R-superlinear convergence of the KKT residues

In this section, we shall derive the asymptotic R-superlinear convergence of the KKT residues on the sequence generated by the ALM for solving (P) under the quadratic growth condition for (D). As explained in section 2.2, this quadratic growth condition is fairly mild for CCCP. Our approach is still within the spirit of Rockafellar's work in [50] in that the ALM is an application of the PPA applied to the dual problem.

Before proceeding, we consider the following three assumptions.

**Assumption 1** *The domain of $h^*$ is an open convex set and $h^*$ is continuously differentiable on* $\mathrm{dom}\, h^*$ *with a globally Lipschitz continuous gradient whose Lipschitz constant is $\lambda_{\nabla h*}$.*

**Assumption 2** *(a) The constraint $\mathcal{B}x \in b + \mathcal{Q}$ takes the form of*

$$\begin{pmatrix} \mathcal{B}_1 \\ \mathcal{B}_2 \end{pmatrix} x = \begin{pmatrix} b_1 \\ b_2 \end{pmatrix} + \begin{pmatrix} \mathcal{Q}_1 \\ \mathcal{Q}_2 \end{pmatrix} \quad \text{(in compatible sizes)},$$

*where $\mathcal{Q}_1 \subseteq \mathbb{Y}_1$ is a polyhedral cone and $\mathcal{Q}_2 \subseteq \mathbb{Y}_2$ is a non-polyhedral cone with nonempty interiors.*
*(b) The function $p$ is split into $p(x) = p_1(x_1) + p_2(x_2)$ for $x = (x_1, x_2) \in \mathbb{X}_1 \times \mathbb{X}_2$, where $p_1 : \mathbb{X}_1 \to (-\infty, +\infty]$ is a proper closed convex polyhedral function and $p_2 : \mathbb{X}_2 \to (-\infty, +\infty]$ is a proper closed convex non-polyhedral function whose domain is a closed convex set with a nonempty interior. Moreover, $p$ is globally Lipschitz continuous on $\mathrm{dom}\, p$ with the Lipschitz constant $\lambda_p$ and $p^*$ is globally Lipschitz continuous on $\mathrm{dom}\, p^*$ with the Lipschitz constant $\lambda_{p*}$.*

**Assumption 3** *The set $\mathrm{SOL}_{\mathrm{D2}}$ is nonempty and the following Robinson constraint qualification (RCQ) of (D2) holds at some $(\bar{w}, \bar{y}, \bar{s}) \in \mathrm{SOL}_{\mathrm{D2}}$ (c.f. [9, Section 3.4.1]):*

$$0 \in \mathrm{int} \left\{ \begin{pmatrix} -c \\ \bar{y} \\ \bar{s} \end{pmatrix} + \begin{pmatrix} \mathcal{A}^* & \mathcal{B}^* & 0 \\ 0 & I_{\mathbb{Y}} & 0 \\ 0 & 0 & I_{\mathbb{X}} \end{pmatrix} \begin{pmatrix} \mathbb{W} \\ \mathbb{Y} \\ \mathbb{X} \end{pmatrix} - \begin{pmatrix} \{0\} \\ \mathcal{Q}^\circ \\ \mathrm{dom}\, p^* \end{pmatrix} \right\},$$

*where $I_{\mathbb{Y}}$ and $I_{\mathbb{X}}$ are the identity operators in $\mathbb{Y}$ and $\mathbb{X}$, respectively.*

Many commonly used functions $h$ and $p$ for CCCP satisfy Assumptions 1 and 2, such as $h$ being any convex quadratic function and $p$ being the indicator function of a closed convex set or any norm function. When $\mathrm{SOL}_{\mathrm{D2}} \neq \emptyset$, Assumption 3 is equivalent to the Mangasarian-Fromovitz constraint qualification in the context of NLP. It is known that under Assumptions 1-3, the optimal solution set $\mathrm{SOL}_{\mathrm{P}}$ to the primal problem is nonempty and bounded [9, Theorem 3.9].

In the rest of this section, we shall derive the asymptotic R-superlinear convergence of the KKT residues for the ALM applied to (P). The following property is useful for developing our main results.

**Lemma 3** *Let $\{(x^k, y^k)\}$ be a sequence generated by the ALM in (1) under criterion (B'). Then for all $k \geq 0$, it holds that*

$$\|y^{k+1} - y^k\| \leq (1 - \eta_k)^{-1} \mathrm{dist}(y^k, \mathrm{SOL}_{\mathrm{D}}).$$

*Proof* By taking $T = T_g$, we know from Proposition 1(b) that for any $\bar{y} \in \mathrm{SOL}_{\mathrm{D}}$,

$$\|y^k - P_k(y^k)\| \leq \|y^k - \bar{y}\|, \quad \forall\, k \geq 0.$$

Hence, $\|y^k - P_k(y^k)\| \leq \text{dist}(y^k, \text{SOL}_D)$ for any $k \geq 0$. Therefore, we have for all $k \geq 0$ that

$$\begin{aligned} \|y^{k+1} - y^k\| &\leq \|y^k - P_k(y^k)\| + \|P_k(y^k) - y^{k+1}\| \\ &\leq \text{dist}(y^k, \text{SOL}_D) + (2\sigma_k(f_k(x^{k+1}) - \inf f_k))^{1/2} \\ &\leq \text{dist}(y^k, \text{SOL}_D) + \eta_k \|y^{k+1} - y^k\|, \end{aligned}$$

where the second term in the second inequality comes from Lemma 1(a), and the third inequality follows from criterion $(B')$. Thus, the conclusion of this lemma follows. $\square$

The following theorem provides the asymptotic Q-superlinear convergence of $\{y^k\}$ and the asymptotic R-superlinear convergence of the primal feasibility, complementarity and primal objective value.

**Theorem 1** *Let $\{(x^k, y^k)\}$ be an infinite sequence generated by the ALM under criterion $(A')$. Then, the sequence $\{y^k\}$ is bounded and converges to some $y^\infty \in \text{SOL}_D$. If criterion $(B')$ is also executed in the ALM and the quadratic growth condition (9) holds at $y^\infty$ with modulus $\kappa_g$, then there exists $\bar{k} \geq 0$ such that for all $k \geq \bar{k}$,*

$$\text{dist}(y^{k+1}, \text{SOL}_D) \leq \mu_k \, \text{dist}(y^k, \text{SOL}_D), \tag{15a}$$

$$\|\Pi_{\mathcal{Q}^\circ}(\mathcal{B}x^{k+1} - b)\| \leq \mu_k' \, \text{dist}(y^k, \text{SOL}_D), \tag{15b}$$

$$|\langle y^{k+1}, \mathcal{B}x^{k+1} - b\rangle| \leq \mu_k'' \, \text{dist}(y^k, \text{SOL}_D), \tag{15c}$$

$$f^0(x^{k+1}) - \inf(P) \leq \mu_k''' \, \text{dist}(y^k, \text{SOL}_D), \tag{15d}$$

*where*

$$\begin{cases} \mu_k := \left[\eta_k + (\eta_k + 1)/\sqrt{1 + \sigma_k^2 \kappa_g^2}\right]/(1 - \eta_k) \to \mu_\infty := 1/\sqrt{1 + \sigma_\infty^2 \kappa_g^2}, \\ \mu_k' := 1/[(1 - \eta_k)\sigma_k] \to \mu_\infty' := 1/\sigma_\infty, \\ \mu_k'' := \|y^{k+1}\|/[(1 - \eta_k)\sigma_k] \to \mu_\infty'' := \|y^\infty\|/\sigma_\infty, \\ \mu_k''' := \left[\eta_k^2 \|y^{k+1} - y^k\| + \|y^{k+1}\| + \|y^k\|\right]/[2(1 - \eta_k)\sigma_k] \to \mu_\infty''' := \|y^\infty\|/\sigma_\infty. \end{cases}$$

*Moreover, $\mu_\infty = \mu_\infty' = \mu_\infty'' = \mu_\infty''' = 0$ if $\sigma_\infty = +\infty$.*

*Proof* The statements on the global convergence just follow from Proposition 2. Next, we prove the results on the rates of convergence.

From Proposition 4 we know that the mapping $T_g^{-1}$ is calm at the origin for $y^\infty$ with modulus $1/\kappa_g$ if the quadratic growth condition (9) holds at $y^\infty$ with modulus $\kappa_g$. Hence, the inequality (15a) is a consequence of Proposition 3(a). It follows from the ALM updating formula $y^{k+1} = \Pi_{\mathcal{Q}^\circ}[y^k + \sigma_k(\mathcal{B}x^{k+1} - b)]$ for all $k \geq 0$ that

$$\mathcal{B}x^{k+1} - b - (1/\sigma_k)(y^{k+1} - y^k) = (1/\sigma_k)\Pi_{\mathcal{Q}}[y^k + \sigma_k(\mathcal{B}x^{k+1} - b)] \in \mathcal{Q},$$

which, implies that for all $k \geq 0$,

$$\begin{aligned} \|\Pi_{\mathcal{Q}^\circ}(\mathcal{B}x^{k+1} - b)\| &= \text{dist}(\mathcal{B}x^{k+1} - b, \mathcal{Q}) \\ &\leq \|\mathcal{B}x^{k+1} - b - [\mathcal{B}x^{k+1} - b - (1/\sigma_k)(y^{k+1} - y^k)]\| = (1/\sigma_k)\|y^{k+1} - y^k\| \end{aligned}$$

and

$$\begin{aligned} |\langle y^{k+1}, \mathcal{B}x^{k+1} - b\rangle| &= |\langle y^{k+1}, (1/\sigma_k)(\Pi_{\mathcal{Q}}[y^k + \sigma_k(\mathcal{B}x^{k+1} - b)] + y^{k+1} - y^k)\rangle| \\ &\leq (1/\sigma_k)\|y^{k+1}\|\|y^{k+1} - y^k\|. \end{aligned}$$

Then the inequalities (15b) and (15c) can be established in view of Lemma 3. Finally, the inequality (15d) follows from (5b) in Proposition 2, criterion $(B')$ and Lemma 3. This completes the proof of the theorem.

Theorem 1 provides fairly general asymptotic superlinear convergence results for the ALM with criteria $(A')$ and $(B')$. Note that, however, it is impractical to execute both $(A')$ and $(B')$ numerically since the values of $\inf f_k$ are not known in general. To circumvent this hurdle, we shall propose verifiable surrogates of criteria $(A')$ and $(B')$ under which we also obtain explicit results on the R-superlinear convergence of the KKT residues for the sequence generated by the ALM.

For any $(w, y, s) \in \mathbb{W} \times \mathbb{Y} \times \mathbb{X}$ and $k \geqslant 0$, denote

$$g_k(y) := g^0(y) - 1/(2\sigma_k)\|y - y^k\|^2, \quad g_k(w, y, s) := g^0(w, y, s) - 1/(2\sigma_k)\|y - y^k\|^2. \quad (16)$$

For any given $k \geqslant 0$ and $y^k \in \mathbb{Y}$, let

$$\begin{cases} \tilde{y}^k(x) := \Pi_{Q^\circ}[y^k + \sigma_k(\mathcal{B}x - b)], \quad \tilde{w}^k(x) := \nabla h(\mathcal{A}x), \\ \tilde{s}^k(x) := \text{Prox}_{p*}[x - (\mathcal{A}^*\tilde{w}^k(x) + \mathcal{B}^*\tilde{y}^k(x) + c)], \\ e^k(x) := x - \text{Prox}_p[x - (\mathcal{A}^*\tilde{w}^k(x) + \mathcal{B}^*\tilde{y}^k(x) + c)], \\ \qquad\quad = \mathcal{A}^*\tilde{w}^k(x) + \mathcal{B}^*\tilde{y}^k(x) + \tilde{s}^k(x) + c \end{cases} \quad x \in \text{dom } f^0. \quad (17)$$

The following lemma provides an upper bound for the duality gap of the ALM subproblem at the $k$-th step. Recall that $f_k(\cdot)$ is the objective function of the ALM subproblem defined in (1).

**Lemma 4** *Let $x \in \text{dom } f^0$. Then for any $k \geqslant 0$,*

$$|f_k(x) - g_k(\tilde{w}^k(x), \tilde{y}^k(x), \tilde{s}^k(x))| \leqslant |\langle x - \tilde{s}^k(x), e^k(x)\rangle| + |p(x) - p(x - e^k(x))|.$$

*Proof* By the definitions of $\tilde{s}^k(x)$ and $e^k(x)$, we have $\tilde{s}^k(x) = \text{Prox}_{p*}[x - e^k(x) + \tilde{s}^k(x)]$, which is also equivalent to $x - e^k(x) \in \partial p^*(\tilde{s}^k(x))$. From [47, Theorems 23.5 & 31.5] and the definition of $\tilde{w}^k(\cdot)$ in (17), we have that

$$h(\mathcal{A}x) + h^*(\tilde{w}^k(x)) = \langle \mathcal{A}x, \tilde{w}^k(x)\rangle, \quad p(x - e^k(x)) + p^*(\tilde{s}^k(x)) = \langle x - e^k(x), \tilde{s}^k(x)\rangle.$$

Hence, we obtain that

$$\begin{aligned}
&|f_k(x) - g_k(\tilde{w}^k(x), \tilde{y}^k(x), \tilde{s}^k(x))| \\
&= |h(\mathcal{A}x) + h^*(\tilde{w}^k(x)) + \langle c, x\rangle + p(x) + p^*(\tilde{s}^k(x)) + \langle \tilde{y}^k(x)/\sigma_k, \tilde{y}^k(x) - y^k + \sigma_k b\rangle| \\
&\leqslant |\langle x, \mathcal{A}^*\tilde{w}^k(x) + c + \tilde{s}^k(x) + \mathcal{B}^*\tilde{y}^k(x)\rangle - \langle e^k(x), \tilde{s}^k(x)\rangle| + |p(x) - p(x - e^k(x))| \\
&\quad + |\langle \tilde{y}^k(x)/\sigma_k, \tilde{y}^k(x) - y^k - \sigma_k(\mathcal{B}x - b)\rangle| \\
&= |\langle x - \tilde{s}^k(x), e^k(x)\rangle| + |p(x) - p(x - e^k(x))| + |\langle \tilde{y}^k(x)/\sigma_k, \Pi_Q[y^k + \sigma_k(\mathcal{B}x - b)]\rangle| \\
&= |\langle x - \tilde{s}^k(x), e^k(x)\rangle| + |p(x) - p(x - e^k(x))|,
\end{aligned}$$

where the last inequality is due to the fact that $\langle \tilde{y}^k(x), \Pi_Q[y^k + \sigma_k(\mathcal{B}x - b)]\rangle = 0$ by the definition of $\tilde{y}^k(\cdot)$. □

The following lemma, which is a direct consequence of [5, Theorem 7], establishes a global upper bound on the distance from a given point to the dual feasible set under the dual Robinson constraint qualification in Assumption 3.

**Lemma 5** *Suppose that Assumptions 1-3 hold. Then there exists a constant $\gamma \geq 1$ such that for any $w \in \operatorname{dom} h^*$, $y \in \mathcal{Q}^\circ$ and $s \in \operatorname{dom} p^*$,*

$$\operatorname{dist}((w,y,s), F_{\mathbf{D}2}) \leq \gamma(1 + \|(w,y,s)\|)\|\mathcal{A}^*w + \mathcal{B}^*y + s + c\|. \tag{18}$$

*Proof* Note that under Assumptions 1-3, one can easily check (c.f. [9, Section 2.3]) that there exists $(\hat{w}, \hat{y}, \hat{s}) \in F_{\mathbf{D}2}$ with $\hat{y} = (\hat{y}_1, \hat{y}_2) \in \mathbb{Y}_1 \times \mathbb{Y}_2$ and $\hat{s} = (\hat{s}_1, \hat{s}_2) \in \mathbb{X}_1 \times \mathbb{X}_2$ such that

$$\hat{y}_2 \in \operatorname{int}(\mathcal{Q}_2^\circ), \quad \hat{s}_2 \in \operatorname{int}(\operatorname{dom} p_2^*).$$

It then follows from [5, Theorem 7] that there exists a positive constant $\bar{\gamma}$ such that for any $w \in \operatorname{dom} h^*$, $y \in \mathcal{Q}^\circ$ and $s \in \operatorname{dom} p^*$,

$$\begin{aligned}
\operatorname{dist}((w,y,s), F_{\mathbf{D}2}) &\leq \bar{\gamma}(1 + \|\Pi_{F_{\mathbf{D}2}}(w,y,s) - (\hat{w},\hat{y},\hat{s})\|)(\|\mathcal{A}^*w + \mathcal{B}^*y + s + c\| + \|s - \Pi_{\operatorname{dom} p^*}(s)\|) \\
&\leq \bar{\gamma}(1 + \|(w,y,s) - (\hat{w},\hat{y},\hat{s})\|)(\|\mathcal{A}^*w + \mathcal{B}^*y + s + c\| + \|s - \Pi_{\operatorname{dom} p^*}(s)\|) \\
&= \bar{\gamma}(1 + \|(w,y,s)\| + \|(\hat{w},\hat{y},\hat{s})\|)\|\mathcal{A}^*w + \mathcal{B}^*y + s + c\|.
\end{aligned}$$

The desired inequality can thus be established with $\gamma = \max\{\bar{\gamma}(1 + \|(\hat{w},\hat{y},\hat{s})\|), 1\}$. $\square$

Equipped with Lemma 5, we are in the position to provide a computable upper bound on the value $f_k(x^{k+1}) - \inf f_k$ used in criteria $(A')$ and $(B')$.

**Proposition 7** *Suppose that Assumptions 1-3 hold. Let $k$ be a fixed but arbitrary nonnegative integer and $\tilde{y}^k(\cdot), \tilde{w}^k(\cdot), \tilde{s}^k(\cdot), e^k(\cdot)$ be the functions defined in (17). Denote $\tilde{z}^k(\cdot) := (\tilde{w}^k(\cdot), \tilde{y}^k(\cdot), \tilde{s}^k(\cdot))$. Suppose that $\{x^{k,j}\}_{j \geq 0}$ is an approximate solution sequence to the augmented Lagrangian subproblem in (1) such that $f_k(x^{k,j}) \to \inf f_k$ and $\|e^k(x^{k,j})\| \to 0$ as $j \to \infty$. Then the sequences $\{x^{k,j}\}_{j \geq 0}$ and $\{\tilde{z}^k(x^{k,j})\}_{j \geq 0}$ are bounded, $\{\tilde{y}^k(x^{k,j})\}_{j \geq 0}$ converges to some point $y^{k,\infty}$ and the following inequality holds:*

$$f_k(x^{k,j}) - g_k(\tilde{z}^k(x^{k,j})) \leq (\|x^{k,j}\| + \|\tilde{s}^k(x^{k,j})\| + \lambda_p)\|e^k(x^{k,j})\|. \tag{19}$$

*Moreover, let $\{t_{k,j}\}_{j \geq 0}$ be a sequence satisfying $1 \geq \sup_{j \geq 0}\{t_{k,j}\} \geq \inf_{j \geq 0}\{t_{k,j}\} > 0$. Then there exists $\bar{j} \geq 0$ such that for any $j \geq \bar{j}$,*

$$\|e^k(x^{k,j})\| \leq \frac{2t_{k,j}}{1 + \|x^{k,j}\| + \|\tilde{z}^k(x^{k,j})\|} \min\left\{\frac{1}{\|\nabla h^*(\tilde{w}^k(x^{k,j}))\| + \|\tilde{y}^k(x^{k,j}) - y^k\|/\sigma_k + 1/\sigma_k}, 1\right\} \tag{20}$$

*and*

$$f_k(x^{k,j}) - \inf f_k \leq \beta^2 t_{k,j}, \tag{21}$$

*where*

$$\beta := \sqrt{2[1 + \lambda_p + \gamma(\|b\| + \lambda_{p*}) + \gamma^2(1 + \lambda_{\nabla h^*})]} \tag{22}$$

*with the constant $\gamma$ given in* (18).

*Proof* From Lemma 1(a) and the assumption that $f_k(x^{k,j}) \to \inf f_k$, we know that the sequence $\{\tilde{y}^j(x^{k,j})\}_{j \geq 0}$ approximately solves $(I + \sigma_k T_g)^{-1}(y^k) = \arg\min\{-g_k(y) \mid y \in \mathbb{Y}\}$, in the sense that

$$(1/2\sigma_k)\|\tilde{y}^k(x^{k,j}) - y^{k,\infty}\|^2 \leq f_k(x^{k,j}) - \inf f_k \to 0 \quad \text{as} \quad j \to \infty,$$

where $y^{k,\infty} \in \mathbb{Y}$ is the unique optimal solution to the strongly convex problem $\min\{-g_k(y) \mid y \in \mathbb{Y}\}$. This further implies that $\{\tilde{y}^k(x^{k,j})\}_{j \geq 0}$ is bounded and converges to $y^{k,\infty}$. By using [49, Theorems 17

& 18], we obtain that $g_k(y^{k,\infty}) = \inf f_k$ and the solution set of the augmented Lagrangian subproblem $\inf f_k$, for which we denote it as $\text{SOL}_{f_k}$, is nonempty and any level set of $f_k$ is bounded. Hence, the sequence $\{x^{k,j}\}_{j\geqslant 0}$ is bounded. The boundedness of $\{\tilde{z}^k(x^{k,j})\}_{j\geqslant 0}$ can be derived from the boundedness of $\{(x^{k,j}, \tilde{y}^k(x^{k,j}))\}_{j\geqslant 0}$ and the local Lipschitz continuity of $\nabla h$.

The inequality (19) is a direct consequence of Lemma 4 and the global Lipschitz continuity of $p$. Since $\|e^k(x^{k,j})\| \to 0$ as $j \to \infty$ and $\inf_{j\geqslant 0}\{t_{k,j}\} > 0$, there exists a nonnegative integer $\bar{j}$ such that (20) holds for any $j \geqslant \bar{j}$. Note that the inequality (20) is valid since $\|e^k(x^{k,j})\| \to 0$ and $\{x^{k,j}\}$, $\{\tilde{z}^k(x^{k,j})\}$ are bounded.

For any $x \in \text{dom}\, f^0$, denote $\bar{z}^k(x) := (\bar{w}^k(x), \bar{y}^k(x), \bar{s}^k(x)) \in \mathbb{W} \times \mathbb{Y} \times \mathbb{X}$ with

$$(\bar{w}^k(x), \bar{y}^k(x), \bar{s}^k(x)) := \Pi_{F_{D2}}(\tilde{z}^k(x)).$$

By the global Lipschitz continuity of $\nabla h^*$ and $p^*$, we obtain that for all $j \geqslant 0$,

$$|g_k(\tilde{z}^k(x^{k,j})) - g_k(\bar{z}^k(x^{k,j}))|$$
$$\leqslant |h^*(\tilde{w}^k(x^{k,j})) - h^*(\bar{w}^k(x^{k,j}))| + |\langle b, \tilde{y}^k(x^{k,j}) - \bar{y}^k(x^{k,j})\rangle| + |p^*(\tilde{s}^k(x^{k,j})) - p^*(\bar{s}^k(x^{k,j}))|$$
$$\quad - (1/2\sigma_k)\big|\|\tilde{y}^k(x^{k,j}) - y^k\|^2 - \|\bar{y}^k(x^{k,j}) - y^k\|^2\big|$$
$$\leqslant |\langle \nabla h^*(\tilde{w}^k(x^{k,j})), \tilde{w}^k(x^{k,j}) - \bar{w}^k(x^{k,j})\rangle| + (\lambda_{\nabla h^*}/2)\|\tilde{w}^k(x^{k,j}) - \bar{w}^k(x^{k,j})\|^2 + \|b\|\|\tilde{y}^k(x^{k,j}) - \bar{y}^k(x^{k,j})\|$$
$$\quad + \lambda_{p^*}\|\tilde{s}^k(x^{k,j}) - \bar{s}^k(x^{k,j})\| + (1/\sigma_k)|\langle \tilde{y}^k(x^{k,j}) - y^k, \tilde{y}^k(x^{k,j}) - \bar{y}^k(x^{k,j})\rangle| + (1/2\sigma_k)\|\tilde{y}^k(x^{k,j}) - \bar{y}^k(x^{k,j})\|^2$$
$$\leqslant \big(\|b\| + \lambda_{p^*} + \|\nabla h^*(\tilde{w}^k(x^{k,j}))\| + (1/\sigma_k)\|\tilde{y}^k(x^{k,j}) - y^k\|\big)\|\tilde{z}^k(x^{k,j}) - \bar{z}^k(x^{k,j})\|$$
$$\quad + \tfrac{1}{2}(\lambda_{\nabla h^*} + 1/\sigma_k)\|\tilde{z}^k(x^{k,j}) - \bar{z}^k(x^{k,j})\|^2.$$

Note that $\|\tilde{z}^k(x^{k,j}) - \bar{z}^k(x^{k,j})\| = \text{dist}(\tilde{z}^k(x^{k,j}), F_{D2})$ and

$$\mathcal{A}^*\tilde{w}^k(x) + \tilde{y}^k(x) + \tilde{s}^k(x) + c = e^k(x), \quad \forall\, x \in \text{dom}\, f^0.$$

By using Lemma 5, we further get

$$|g_k(\tilde{z}^k(x^{k,j})) - g_k(\bar{z}^k(x^{k,j}))|$$
$$\leqslant \gamma\big(\|b\| + \lambda_{p^*} + \|\nabla h^*(\tilde{w}^k(x^{k,j}))\| + (1/\sigma_k)\|\tilde{y}^k(x^{k,j}) - y^k\|\big)(1 + \|\tilde{z}^k(x^{k,j})\|)\|e^k(x^{k,j})\|$$
$$\quad + \tfrac{1}{2}\gamma^2(\lambda_{\nabla h^*} + 1/\sigma_k)(1 + \|\tilde{z}^k(x^{k,j})\|)^2\|e^k(x^{k,j})\|^2.$$

The above inequality, together with (19), implies that for $j \geqslant \bar{j}$,

$$f_k(x^{k,j}) - \inf f_k = f_k(x^{k,j}) - \sup_y\{g_k(w,y,s) \mid (w,y,s) \in F_{D2}\}$$
$$\leqslant f_k(x^{k,j}) - g_k(\bar{z}^k(x^{k,j})) = f_k(x^{k,j}) - g_k(\tilde{z}^k(x^{k,j})) + g_k(\tilde{z}^k(x^{k,j})) - g_k(\bar{z}^k(x^{k,j}))$$
$$\leqslant \big(1 + \lambda_p + \gamma(\|b\| + \lambda_{p^*}) + \gamma^2 \lambda_{\nabla h^*}\big)(1 + \|x^{k,j}\| + \|\tilde{z}^k(x^{k,j})\|)\|e^k(x^{k,j})\|$$
$$\quad + \big(\gamma\|\nabla h^*(\tilde{w}^k(x^{k,j}))\| + \gamma\|\tilde{y}^k(x^{k,j}) - y^k\|/\sigma_k + \gamma^2/\sigma_k\big)(1 + \|x^{k,j}\| + \|\tilde{z}^k(x^{k,j})\|)\|e^k(x^{k,j})\|,$$

where in deriving the second inequality, we used the fact in (20) that $(1+\|x^{k,j}\|+\|\tilde{z}^k(x^{k,j})\|)\|e^k(x^{k,j})\| \leqslant 2$. Now by using the inequality (20) and the fact that $\gamma \geqslant 1$, we can show that (21) holds for any $j \geqslant \bar{j}$. This completes the proof.

Let $\{\hat{\varepsilon}_k\}$ and $\{\hat{\eta}_k\}$ be two given positive summable sequences. Suppose that for some $k \geq 0$, $y^k \in \mathcal{Q}^\circ$ is not an optimal solution to (D). Let $\{x^{k,j}\}_{j \geq 0}$ be a sequence satisfying $f_k(x^{k,j}) \to \inf f_k$ and $\|e^k(x^{k,j})\| \to 0$ as $j \to \infty$. Then the limit $y^{k,\infty}$ of $\{\tilde{y}^k(x^{k,j})\}_{j \geq 0}$ must be different from $y^k$ because otherwise $y^k$ and any accumulation point of $\{x^{k,j}\}_{j \geq 0}$ would form a KKT solution point to (P), contradicting the assumption that $y^k$ does not solve (D). This observation allows us to apply Proposition 7 with $\{t_{k,j}\}$ being either $\{\hat{\varepsilon}_k^2/2\sigma_k\}$ or $\{\hat{\eta}_k^2/(2\sigma_k)\|\tilde{y}^k(x^{k,j}) - y^k\|^2\}$, both of which are bounded away from 0 as $y^k$ does not solve (D). Thus, we can use the following two criteria to replace $(A')$ and $(B')$, respectively:

$$(A'') \quad \|e^{k+1}\| \leq \frac{\hat{\varepsilon}_k^2/\sigma_k}{1 + \|x^{k+1}\| + \|z^{k+1}\|} \min\left\{\frac{1}{\|\nabla h^*(w^{k+1})\| + \|y^{k+1} - y^k\|/\sigma_k + 1/\sigma_k}, 1\right\},$$

$$(B'') \quad \|e^{k+1}\| \leq \frac{(\hat{\eta}_k^2/\sigma_k)\|y^{k+1} - y^k\|^2}{1 + \|x^{k+1}\| + \|z^{k+1}\|} \min\left\{\frac{1}{\|\nabla h^*(w^{k+1})\| + \|y^{k+1} - y^k\|/\sigma_k + 1/\sigma_k}, 1\right\},$$

where

$$w^{k+1} := \tilde{w}(x^{k+1}), \ y^{k+1} := \tilde{y}(x^{k+1}), \ s^{k+1} := \tilde{s}(x^{k+1}), \ z^{k+1} := (w^{k+1}, y^{k+1}, s^{k+1}), \ e^{k+1} := e^k(x^{k+1}).$$

Intuitively, for the augmented Lagrangian subproblem $\inf f_k$, its optimality condition can be written as

$$x - \text{Prox}_p\left[x - (\mathcal{A}^*\nabla h(\mathcal{A}x) + c + \mathcal{B}^*\Pi_{\mathcal{Q}^\circ}[y^k - \sigma_k(b - \mathcal{B}x)])\right] = 0, \quad x \in \text{dom } f^0.$$

Thus, if the $(k+1)$-th subproblem in the ALM is solved exactly, one must have $e^{k+1} = 0$, indicating that both $(A'')$ and $(B'')$ hold automatically. In fact, criteria $(A'')$ and $(B'')$ essentially require the residue of the above nonsmooth equation at the current iteration to be sufficiently small, whereas the original criteria $(A')$ and $(B')$ proposed by Rockafellar in [50] ask the gap between the current objective function value and the optimal objective function value to be sufficiently small. Proposition 7 says that under Assumptions 1-3, criteria $(A')$ and $(B')$ are satisfied with $\varepsilon_k = \beta\hat{\varepsilon}_k$ and $\eta_k = \beta\hat{\eta}_k$ as long as $(A'')$ and $(B'')$ are true. As far as we know, these easy-to-implement criteria for the ALM have never been discovered before. It is also worth mentioning that the computation of $e^{k+1}$ requires not much extra costs in numerical implementations, as shall be demonstrated in the next section.

Denote the natural map associated with the KKT optimality condition (7) as

$$R^{\text{nat}}(x, y) := \begin{pmatrix} x - \text{Prox}_p[x - (\mathcal{A}^*\nabla h(\mathcal{A}x) + \mathcal{B}^*y + c)] \\ y - \Pi_{\mathcal{Q}^\circ}[y - (b - \mathcal{B}x)] \end{pmatrix}, \quad (x, y) \in \text{dom } f^0 \times \mathbb{Y}. \tag{23}$$

It can be easily seen that $x \in \text{SOL}_\text{P}$ and $y \in \text{SOL}_\text{D}$ if and only if $R^{\text{nat}}(x, y) = 0$. The following theorem establishes the much promised global convergence and the asymptotic R-superlinear convergence of the KKT residues in terms of $\|R^{\text{nat}}(x^k, y^k)\|$ under criteria $(A'')$ and $(B'')$.

**Theorem 2** *Suppose that Assumptions 1-3 hold. Let $\{(x^k, y^k)\}$ be an infinite sequence generated by the ALM in (1) under criterion $(A'')$. Then the sequence $\{y^k\}$ is bounded and converges to some $y^\infty \in \text{SOL}_\text{D}$, and the sequence $\{x^k\}$ is also bounded with all of its accumulation points in $\text{SOL}_\text{P}$.*

*If criterion $(B'')$ is also executed in the ALM and the quadratic growth condition (9) holds at $y^\infty$ with modulus $\kappa_g$, then there exists $k' \geq 0$ such that for all $k \geq k'$, $\beta\hat{\eta}_k < 1$ and*

$$\text{dist}(y^{k+1}, \text{SOL}_\text{D}) \leq \theta_k \text{dist}(y^k, \text{SOL}_\text{D}), \tag{24a}$$

$$\|R^{\text{nat}}(x^{k+1}, y^{k+1})\| \leq \theta_k' \text{dist}(y^k, \text{SOL}_\text{D}), \tag{24b}$$

*where*

$$\begin{cases} \theta_k := \left[\beta\hat{\eta}_k + (\beta\hat{\eta}_k + 1)/\sqrt{1 + \sigma_k^2\kappa_g^2}\right]/(1 - \beta\hat{\eta}_k) \to \theta_\infty := 1/\sqrt{1 + \sigma_\infty^2\kappa_g^2}, \\ \theta_k' := [\max\{1, 1/\sigma_k\} + (\hat{\eta}_k^2/\sigma_k)\|y^{k+1} - y^k\|]/(1 - \beta\hat{\eta}_k) \to \theta_\infty' := \max\{1, 1/\sigma_\infty\}, \end{cases}$$

*where the constant $\beta$ is given in* (22). *Moreover, when $\mathcal{Q} = \{0\}$, i.e., the constraint in* (P) *is $\mathcal{B}x = b$, the term $\theta_k'$ in* (24b) *can be replaced by*

$$\theta_k' := [1/\sigma_k + (\hat{\eta}_k^2/\sigma_k)\|y^{k+1} - y^k\|]/(1 - \beta\hat{\eta}_k) \to \theta_\infty' := 1/\sigma_\infty.$$

*Proof* All the statements except the inequality (24b) follow from Theorem 1 and the inequality (21) in Proposition 7. Noting that under criterion $(B'')$, we have

$$\|x^{k+1} - \text{Prox}_p[x^{k+1} - (\mathcal{A}^*\nabla h(\mathcal{A}x^{k+1}) + \mathcal{B}^*y^{k+1} + c)]\| = \|e^{k+1}\| \leqslant (\hat{\eta}_k^2/\sigma_k)\|y^{k+1} - y^k\|^2, \quad \forall\, k \geqslant 0. \quad (25)$$

By using [58, Lemma 2] and (5a) in Proposition 2, we have for any $\sigma_k \geqslant 1$ that

$$\begin{aligned} & \|y^{k+1} - \Pi_{\mathcal{Q}^\circ}[y^{k+1} - (b - \mathcal{B}x^{k+1})]\| \\ & \leqslant \|y^{k+1} - \Pi_{\mathcal{Q}^\circ}[y^{k+1} - \sigma_k(b - \mathcal{B}x^{k+1})]\| \\ & = \|\Pi_{\mathcal{Q}^\circ}[y^k + \sigma_k(\mathcal{B}x^{k+1} - b)] - \Pi_{\mathcal{Q}^\circ}[y^{k+1} + \sigma_k(\mathcal{B}x^{k+1} - b)]\| \\ & \leqslant \|y^{k+1} - y^k\| \end{aligned} \quad (26)$$

and by using [24] and (5a) in Proposition 2, we have for any $0 < \sigma_k < 1$ that

$$\begin{aligned} & \|y^{k+1} - \Pi_{\mathcal{Q}^\circ}[y^{k+1} - (b - \mathcal{B}x^{k+1})]\| \\ & \leqslant (1/\sigma_k)\|y^{k+1} - \Pi_{\mathcal{Q}^\circ}[y^{k+1} - \sigma_k(b - \mathcal{B}x^{k+1})]\| \\ & = (1/\sigma_k)\|\Pi_{\mathcal{Q}^\circ}[y^k + \sigma_k(\mathcal{B}x^{k+1} - b)] - \Pi_{\mathcal{Q}^\circ}[y^{k+1} + \sigma_k(\mathcal{B}x^{k+1} - b)]\| \\ & \leqslant (1/\sigma_k)\|y^{k+1} - y^k\|. \end{aligned} \quad (27)$$

Then, from (26) and (27), we obtain that for any $k \geqslant 0$

$$\|y^{k+1} - \Pi_{\mathcal{Q}^\circ}[y^{k+1} - (b - \mathcal{B}x^{k+1})]\| \leqslant \max\{1, 1/\sigma_k\}\|y^{k+1} - y^k\|.$$

Thus, by using (25) and Lemma 3, we know that there exists $k' \geqslant 0$ such that for all $k \geqslant k'$, $\beta\hat{\eta}_k < 1$ and (24b) holds.

Finally, by observing that for the equality constrained case with $\mathcal{Q} = \{0\}$ in (P), it holds that for all $k \geqslant 0$,

$$\|y^{k+1} - \Pi_{\mathcal{Q}^\circ}[y^{k+1} - (b - \mathcal{B}x^{k+1})]\| = \|b - \mathcal{B}x^{k+1}\| = (1/\sigma_k)\|y^{k+1} - y^k\|,$$

we can complete the proof of this theorem. □

Below we make a couple of remarks regarding the convergence rates proven in Theorems 1 and 2 for the ALM applied to CCCP.

*Remark 1* Under the dual quadratic growth condition at $y^\infty$, the asymptotic Q-(super)linear convergence of the dual sequence $\{y^k\}$ in Theorem 1 is an extension of the results established in [50,37], while the asymptotic R-(super)linear convergence rates of the primal feasibility, complementarity and primal objective value in Theorem 1, and of the KKT residues in Theorem 2, to the best of our knowledge, are only proven here. The latter results resolve the long outstanding open question of the ALM on the lack of convergence rates for the KKT residues when the KKT solution mapping does not possess the restrictive calmness condition. They also reveal that the ALM is indeed a dual-type method, in the sense that its KKT convergence rates can be derived by solely using the property of the dual solution mapping. This feature distinguishes the ALM from the primal-dual-type methods, such as the alternating direction method of multipliers (ADMM), for which the convergence rate is derived under proper assumptions on the KKT solution mapping, see, e.g., [26]. As mentioned in section 3, the calmness of the KKT solution mapping, which is a strictly stronger condition than that of the dual solution mapping, is difficult to be satisfied for non-polyhedral CCCP. Hence, the above arguments partially explain why the ADMM does not perform well for many challenging CCCP problems, as shown by the extensive numerical results conducted in [59,36,14].

*Remark 2* All the convergence rates proven in Theorems 1 and 2 become asymptotically superlinear if the penalty parameters $\sigma_k \to +\infty$. However, in numerical computations one never needs $\sigma_k$ converging to $+\infty$. Instead one can just progressively choose $\sigma_k$ to be large enough, such as $\sigma_k \approx 1/\kappa_g$ with $\kappa_g$ being the dual quadratic growth modulus, to achieve a fast linear rate. Of course, in general one does not know $\kappa_g$ in practice, and hence the adjustment of $\sigma_k$ to achieve fast convergence is always an important issue in the practical implementation of the ALM.

## 5 Applications to linear and convex quadratic SDP

In this section, we will illustrate the usefulness of the results obtained in the last section in the context of linear and convex quadratic SDP problems. With the rich structure exhibited in this class of problems, we are able to gain more insights on the superior properties of the ALM. The convex quadratic SDP problem and its dual take the following forms:

$$\min_{x=(x_1,x_2,x_3)} f^0(x) := \frac{1}{2}\langle x_1, \mathcal{H}x_1\rangle - \langle b, x_2\rangle + \delta_{\mathbb{S}^n_+}(x_3) \quad \text{(QSDP-P)}$$
$$\text{s.t.} \quad -\mathcal{H}x_1 + \mathcal{E}^* x_2 + x_3 = C, \ x_1 \in \text{Ran}(\mathcal{H})$$

and

$$\max_X g^0(X) := -\frac{1}{2}\langle X, \mathcal{H}X\rangle - \langle C, X\rangle \quad \text{(QSDP-D)}$$
$$\text{s.t.} \quad \mathcal{E}X = b, \quad X \in \mathbb{S}^n_+,$$

where $\mathcal{H} : \mathbb{S}^n \to \mathbb{S}^n$ is a self-adjoint positive semidefinite linear operator, $\mathcal{E} : \mathbb{S}^n \to \mathbb{R}^m$ is a linear map, $C \in \mathbb{S}^n$ and $b \in \mathbb{R}^m$ are given data, and $\text{Ran}(\mathcal{H})$ is the range space of $\mathcal{H}$. Here, we swap the primal and dual formats from the conventional treatments such as in [36] to make our discussions consistent with the previous sections. Problem (QSDP-P) is a special instance of (P) by taking the functions $h$ and $p$ as

$$h(x) = \frac{1}{2}\langle x_1, \mathcal{H}x_1\rangle, \quad p(x) = \delta_{\mathbb{S}^n_+}(x_3), \quad x = (x_1, x_2, x_3) \in \mathbb{X} := \text{Ran}(\mathcal{H}) \times \mathbb{R}^m \times \mathbb{S}^n.$$

In the solvers SDPNAL/SDPNAL+ (for solving linear SDP with $\mathcal{H} = 0$) [60,59] and QSDPNAL [36], the authors designed ALMs for solving (QSDP-P) that employed the semismooth Newton-CG method to obtain high quality approximate solutions for the inner subproblems. Specifically, given an initial point $X^0 \in \mathbb{S}^n_+$ and a positive scalar sequence $\sigma_k \uparrow \sigma_\infty \leq +\infty$, the $k$-th iteration of the ALM is given by

$$\begin{cases} x^{k+1} \approx \arg\min_{x \in \mathbb{X}} f_k(x) := f^0(x) + \frac{1}{2\sigma_k}(\|X^k + \sigma_k(-\mathcal{H}x_1 + \mathcal{E}^*x_2 + x_3 - C)\|^2 - \|X^k\|^2), \\ X^{k+1} = X^k + \sigma_k(-\mathcal{H}x_1 + \mathcal{E}^*x_2 + x_3 - C). \end{cases}$$

We have the following sufficient conditions of the dual quadratic growth condition for the convex quadratic SDP problems, which can be obtained from [14, Theorems 13 & 18].

**Proposition 8** *Assume that the KKT system to* (QSDP-P) *has at least one solution. Then, the quadratic growth condition for* (QSDP-D) *holds at a dual optimal solution* $\overline{X} \in \mathbb{S}^n$ *if either there exists a KKT point* $(\hat{x}, \hat{X}) \in \mathbb{X} \times \mathbb{S}^n$ *such that* $\text{rank}(\hat{X}) + \text{rank}(\hat{x}_3) = n$ *($\hat{X}$ may be different from $\overline{X}$), or the following SOSC holds at $\overline{X}$:*

$$\sup_{(\bar{x}_1, \bar{x}_2, \bar{x}_3) \in \text{SOL}_P} \langle D, \mathcal{H}D \rangle + 2\langle \bar{x}_3, D\overline{X}^\dagger D \rangle > 0, \quad \forall D \in \mathcal{C}(\overline{X}) \setminus \{0\},$$

*where $\overline{X}^\dagger$ is the Moore-Penrose pseudoinverse of $\overline{X}$ and $\mathcal{C}(\overline{X})$ is the critical cone of the dual problem at $\overline{X}$, i.e.,*

$$\mathcal{C}(\overline{X}) = \left\{ D \in \mathbb{S}^n \mid \mathcal{E}D = 0, \ D \in \mathcal{T}_{\mathbb{S}^n_+}(\overline{X}), \ \langle D, \mathcal{H}\overline{X} + C \rangle = 0 \right\}.$$

Proposition 8 says that for the convex quadratic SDP problem, if a partial strict complementarity solution associated with the non-polyhedral constraint $X \in \mathbb{S}^n_+$ exists, then the dual quadratic growth condition holds at any dual optimal solution. In fact, we have a relatively complete picture on the relationships between the following concepts and properties associated with convex quadratic SDP: the Lipschitz continuity (Lip.) and robust isolated calmness (r.iso.calm) of the primal/dual/KKT solution mappings at the origin, the quadratic growth condition (q-grow.) for the primal/dual problems at optimal solution points, the SOSC for the primal/dual problems, the extended strict Robinson constraint qualification (eSRCQ) for the primal/dual problems (for its definition, see [27, Definition 5.1]) and the existence of a partial strict complementarity KKT solution ($\exists$ strict comp.sol.). They are summarized in the diagram below.

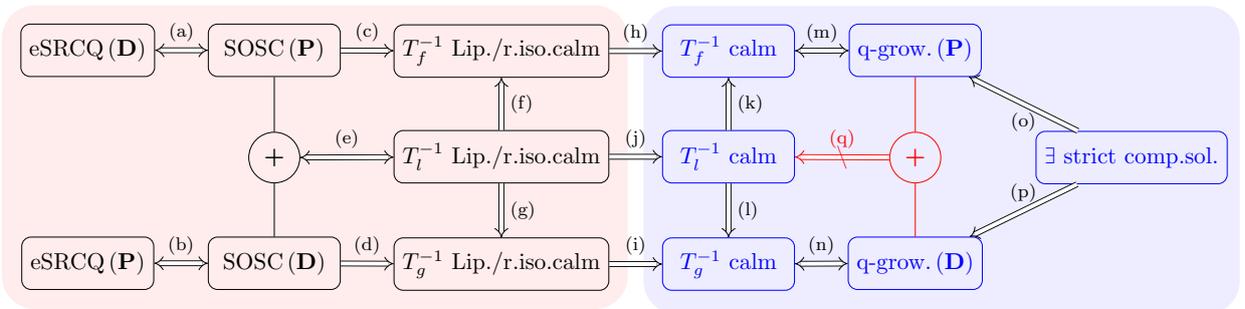

**Fig. 2** A diagram of the Lipschitzian-like properties for convex quadratic SDP (the first three columns refer to the cases with a unique solution while the last three columns refer to the cases with possibly multiple solutions.)

In the above diagram, the relations $(a)$ and $(b)$ are from [27, Propositions 5.3 & 5.4]; the implications $(c)$, $(d)$, $(o)$ and $(p)$ are results of [14, Theorems 12 & 13] (see also Proposition 8 in the above); the equivalence in $(e)$ between the SOSCs for both (P) and (D) and the Lipschitz continuity or the robust isolated calmness of $T_l^{-1}$ at the origin is given by Proposition 6 and [27, Theorem 5.1]; the implications $(f)$-$(l)$ can be directly obtained by definitions; the relations $(m)$ and $(n)$ are stated in Proposition 4 and the negated implication $(q)$ is demonstrated by Example 2.

Figure 2 further explains that the dual quadratic growth condition is fairly mild, much weaker than the calmness of $T_l^{-1}$ at the origin. As suggested by Theorem 2, the convergence rate of the KKT residues for $\{(x^k, X^k)\}$ is asymptotically R-superlinear under the quadratic growth condition for (QSDP-D), while the weakest known assumption to ensure the linear convergence rate of the ADMM for solving (QSDP-P) is the calmness of the KKT solution mapping $T_l^{-1}$ at the origin for a KKT point [26].

In the following, we shall illustrate how to implement criteria $(A'')$ and $(B'')$ if the subproblems in the ALM are solved by the semismooth Newton-CG method without incurring significant extra computational costs. For any $(x_1, x_2) \in \text{Ran}(\mathcal{H}) \times \mathbb{R}^m$ and $k \geq 0$, denote

$$\psi_k(x_1, x_2) := \frac{1}{2}\langle x_1, \mathcal{H}x_1\rangle - \langle b, x_2\rangle + \frac{1}{2\sigma_k}(\|\Pi_{\mathbb{S}^n_+}[X^k + \sigma_k(-\mathcal{H}x_1 + \mathcal{E}^*x_2 - C)]\|^2 - \|X^k\|^2).$$

One can easily check that if $(\tilde{x}_1, \tilde{x}_2, \tilde{x}_3) \in \arg\min\{f_k(x) \mid x \in \mathbb{X}\}$, we have

$$\begin{cases} (\tilde{x}_1, \tilde{x}_2) \in \arg\min\{\psi_k(x_1, x_2) \mid (x_1, x_2) \in \text{Ran}(\mathcal{H}) \times \mathbb{R}^m\}, \\ \tilde{x}_3 = (1/\sigma_k)\Pi_{\mathbb{S}^n_+}[X^k + \sigma_k(-\mathcal{H}\tilde{x}_1 + \mathcal{E}^*\tilde{x}_2 - C)]. \end{cases} \quad (28)$$

To solve the above optimization problem associated with $(x_1, x_2)$ inexactly, it suffices to approximately solve the following nonsmooth equation by the semismooth Newton-CG method:

$$\nabla \psi_k(x_1, x_2) = \begin{pmatrix} \mathcal{H}x_1 - \mathcal{H}\Pi_{\mathbb{S}^n_+}[X^k + \sigma_k(-\mathcal{H}x_1 + \mathcal{E}^*x_2 - C)] \\ -b + \mathcal{E}\Pi_{\mathbb{S}^n_+}[X^k + \sigma_k(-\mathcal{H}x_1 + \mathcal{E}^*x_2 - C)] \end{pmatrix} = 0, \quad (29)$$

where $(x_1, x_2) \in \text{Ran}(\mathcal{H}) \times \mathbb{R}^m$. After obtaining an approximate solution $(x_1^{k+1}, x_2^{k+1}) \in \text{Ran}(\mathcal{H}) \times \mathbb{R}^m$ to the above equation, we set

$$x_3^{k+1} = (1/\sigma_k)\Pi_{\mathbb{S}^n_+}[X^k + \sigma_k(-\mathcal{H}x_1^{k+1} + \mathcal{E}^*x_2^{k+1} - C)].$$

Direct computations show that the vector $e^{k+1}$ used in criteria $(A'')$ and $(B'')$ is given by

$$e^{k+1} = \begin{pmatrix} \nabla\psi_k(x_1^{k+1}, x_2^{k+1}) \\ 0 \end{pmatrix} \quad \text{with} \quad x^{k+1} = (x_1^{k+1}, x_2^{k+1}, x_3^{k+1}).$$

If a semismooth Newton-CG method is adopted to solve the equation (29), which is exactly the method employed in the solvers SDPNAL [60], SDPNAL+ [59] and QSDPNAL [36], the stopping criteria $(A'')$ and $(B'')$ then turn out to be

$$(\widetilde{A}'') \quad \|\nabla\psi_k(x_1^{k+1}, x_2^{k+1})\| \leq \frac{\tilde{\varepsilon}_k^2/\sigma_k}{1 + \|x^{k+1}\| + \|X^{k+1}\|} \min\left\{\frac{1}{\|\mathcal{H}X^{k+1}\| + \|X^{k+1} - X^k\|/\sigma_k + 1/\sigma_k}, 1\right\},$$

$$(\widetilde{B}'') \quad \|\nabla\psi_k(x_1^{k+1}, x_2^{k+1})\| \leq \frac{(\tilde{\eta}_k^2/\sigma_k)\|X^{k+1} - X^k\|^2}{1 + \|x^{k+1}\| + \|X^{k+1}\|} \min\left\{\frac{1}{\|\mathcal{H}X^{k+1}\| + \|X^{k+1} - X^k\|/\sigma_k + 1/\sigma_k}, 1\right\}$$

with the given positive summable sequences $\{\widetilde{\varepsilon}_k\}$ and $\{\widetilde{\eta}_k\}$. Since the sequence $\{x^k\}$ generated by the semismooth Newton-CG method satisfies $f_k(x^{k+1}) \to \inf f_k$, under the RCQ for (QSDP-D), criteria $(\widetilde{A}'')$ and $(\widetilde{B}'')$ are achievable due to Proposition 7. Interestingly, to check criteria $(\widetilde{A}'')$ and $(\widetilde{B}'')$ do not need much extra computational costs, since the value $\nabla \psi_k(x_1^{k+1}, x_2^{k+1})$ shall be used in the next iteration step of the inner semismooth Newton-CG method, unless the point $(x_1^{k+1}, x_2^{k+1})$ is already accepted as an approximate solution to the subproblem in (28).

The natural map associated with the KKT optimality condition for (QSDP-P), in the sense of (23), takes the form of

$$R(x, X) := \begin{pmatrix} \mathcal{H}x_1 - \mathcal{H}X \\ -b + \mathcal{E}X \\ x_3 - \Pi_{\mathbb{S}_+^n}(x_3 - X) \\ \mathcal{H}x_1 - \mathcal{E}^* x_2 - x_3 + C \end{pmatrix}, \quad x = (x_1, x_2, x_3) \in \text{Ran}(\mathcal{H}) \times \mathbb{R}^m \times \mathbb{S}_+^n, \ X \in \mathbb{S}^n.$$

Based on Theorem 2, under the Robinson constraint qualification and the mild dual quadratic growth condition and criteria $(\widetilde{A}'')$ and $(\widetilde{B}'')$, we know that the sequence $\{X^k\}$ generated by the ALM converges asymptotically Q-superlinearly, and the KKT residues $\{\|R(x^k, X^k)\|\}$ converges asymptotically R-superlinearly. These much desired fast convergence rates have been observed in the solvers SDPNAL [60], SDPNAL+ [59] and QSDPNAL [36] for solving linear and convex quadratic SDP problems and actually can be used to explain the superior performance of these solvers.

In the last part of this section, we shall conduct some numerical experiments to verify the derived convergence rates. Consider the following weighted nearest correlation matrix problem:

$$\begin{aligned} \min_{X} \ & \tfrac{1}{2} \|H \circ (X - G)\|^2 \\ \text{s.t.} \ & X_{ii} = 1, \ i = 1, 2, \ldots, n, \quad X \in \mathbb{S}_+^n, \end{aligned} \tag{30}$$

where $G \in \mathbb{S}^n$ is an observed sample correlation matrix, $H \in \mathbb{S}^n$ is a given nonnegative weight matrix and $\circ$ denotes the Hadamard product, i.e., $(A \circ B)_{ij} = A_{ij} B_{ij}$ for any $A, B \in \mathbb{S}^n$. This model has been widely used in finance for estimating sample correlation matrices with missing data, where a typical choice of the weight matrix $H$ is to ask $H_{ij} = 1$ if $G_{ij}$ is observed and $H_{ij} = 0$ if $G_{ij}$ is missing [29]. Other examples of the matrix $H$ in finance can be found in [8].

Define $\mathcal{H}(W) = H \circ H \circ W$ for $W \in \mathbb{S}^n$, and $C = H \circ H \circ G$. The dual of (30) can be written as

$$\begin{aligned} \min_{y, W, S} \ & \tfrac{1}{2} \langle W, \mathcal{H}W \rangle - \sum_{i=1}^n y_i + \delta_{\mathbb{S}_+^n}(S), \\ \text{s.t.} \ & \text{diag}(y) - \mathcal{H}W + S = C, \ W \in \text{Ran}(\mathcal{H}), \end{aligned} \tag{31}$$

where $\text{diag}(y)$ is the diagonal matrix with the vector $y$ as its diagonal. Obviously (31) is in the form of (QSDP-P). We can apply the ALM in (28) to solve (31) with the subproblems being solved by the semismooth Newton-CG method. One nice feature of the nearest correlation matrix problem is that the constraints in (30) are always nondegenerate, making the semismooth Newton-CG method converging at least Q-superlinearly [43].

In our numerical experiments, we take the matrix $G$ to be the indefinite symmetric matrices constructed from stock data by the investment company Orbis [30] (one with the matrix dimensions $1399 \times 1399$ and the other with dimensions $3210 \times 3210$), which are available at https://github.com/higham/matrices-correlation-invalid. We randomly set some entries $H_{ij} = 0$ (the corresponding $G_{ij}$ are thus treated as "missing" from the the observations) and the other entries $H_{ij} = 1$. Figure 3 shows the semi-log curves of the KKT residues versus the ALM iteration count, where one can easily observe the R-superlinear convergence.

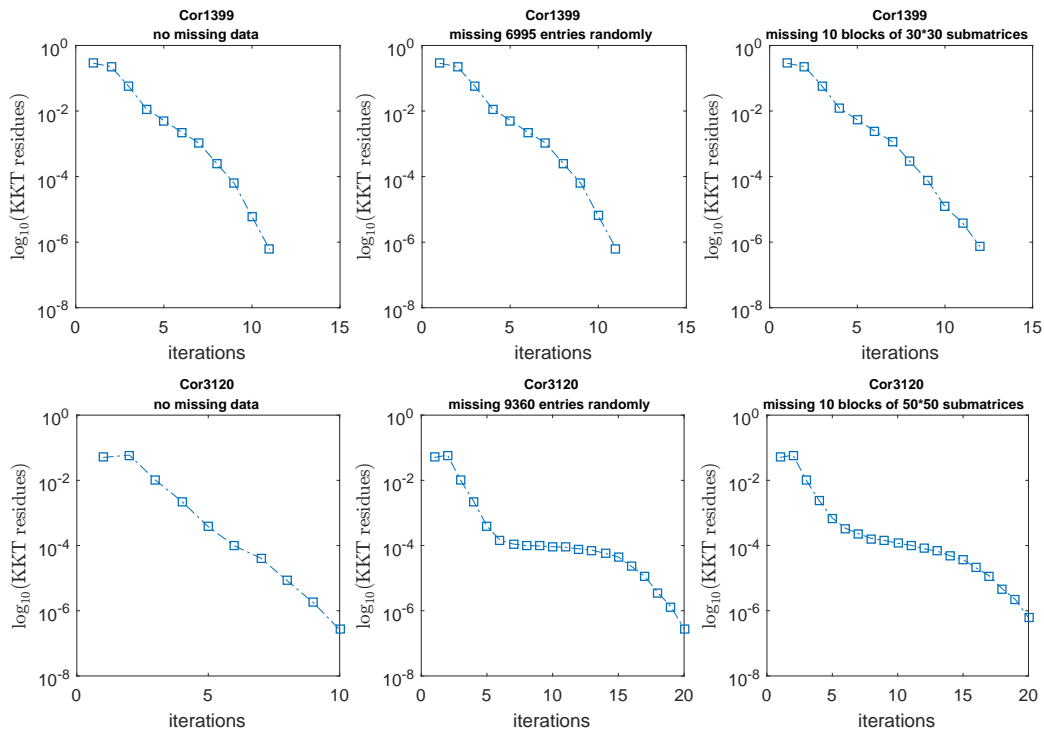

**Fig. 3** The KKT residues of the ALM for solving the dual problem (31) of the H-weighted nearest correlation matrix problem (30) with missing data.

## 6 Concluding discussions

In this paper, we have established the asymptotic R-superlinear convergence of the KKT residues for the iterates generated by the ALM for solving CCCP, under a fairly mild quadratic growth condition on the dual problem, for which neither the primal nor the dual solution is required to be unique. The obtained result fully explains the numerical success of the solvers SDPNAL [60], SDPNAL+ [59] and QSDPNAL [36] for solving linear and convex quadratic SDP problems. We believe that the research presented in this paper has provided a practical guide for designing efficient general large scale CCCP solvers in the future. One question that has not been answered here is whether the primal sequence generated by the ALM for solving CCCP can also converge asymptotically superlinearly without the upper Lipschitz continuity of the KKT solution mapping at the origin. Though this question is mainly of theoretical importance given what have been achieved in this paper on the fast convergence rates of the KKT residues, its answer may help deepen our understanding on the ALM.

## References

1. Alizadeh, F., Haeberly, J.-P.A., Overton, M.L.: Complementarity and nondegeneracy in semidefinite programming. Math. Program. **77**(2), 111–128 (1997)
2. Alves, M.M., Svaiter, B.F.: A note on Fejér-monotone sequences in product spaces and its applications to the dual convergence of augmented Lagrangian methods. Math. Program. **155**(1), 613–616 (2016)
3. Artacho, F.J.A., Geoffroy, M.H.: Characterization of metric regularity of subdifferentials. J. Convex Anal. **15**(2), 365–380 (2008)

4. Attouch, H., Soueycatt, M.: Augmented Lagrangian and proximal alternating direction methods of multipliers in Hilbert spaces. Applications to games, PDE's and control. Pac. J. Optim. **5**(1), 17–37 (2009)
5. Bauschke, H.H., Borwein, J.M., Li, W.: Strong conical hull intersection property, bounded linear regularity, Jameson's property (G), and error bounds in convex optimization. Math. Program. **86**(1), 135–160 (1999)
6. Bergounioux, M.: Use of augmented Lagrangian methods for the optimal control of obstacle problems. J. Optim. Theory Appl. **95**(1), 101–126 (1997)
7. Bertsekas, D.: Constrained Optimization and Lagrange Multipliers Methods. Academic Press, New York (1982)
8. Bhansali, V., Wise, B.: Forecasting portfolio risk in normal and stressed market. J. Risk, **4**(1), 91–106 (2001)
9. Bonnans, J.F., Shapiro, A.: Perturbation Analysis of Optimization Problems. Springer, New York (2000)
10. Chen, C.H., Liu, Y.J., Sun, D.F., Toh, K.-C.: A semismooth Newton-CG dual proximal point algorithm for matrix spectral norm approximation problems. Math. Program. **155**(1), 435–470 (2016)
11. Conn, A.R., Gould, N., Sartenaer, A., Toint, P.L.: Convergence Properties of an Augmented Lagrangian Algorithm for Optimization with a Combination of General Equality and Linear Constraints. SIAM J. Optim. **6**(3), 674–703 (1996)
12. Conn, A.R., Gould, N., Toint, P.L.: A globally convergent augmented Lagrangian algorithm for optimization with general constraints and simple bounds. SIAM J. Numer. Anal. **28**(2), 545–572 (1991)
13. Contesse-Becker, L.: Extended convergence results for the method of multipliers for non-strictly binding inequality constraints. J. Optim. Theory Appl. **79**(2), 273-310 (1993)
14. Cui, Y., Ding, C., Zhao, X.Y.: Quadratic growth conditions for convex matrix optimization problems associated with spectral functions. arXiv:1702.03262, to appear in SIAM J. Optim. (2017)
15. Cui, Y., Sun, D.F.: A complete characterization on the robust isolated calmness of the nuclear norm regularized convex optimization problems. arXiv:1702.05914 (2017)
16. Ding, C.: An Introduction to a Class of Matrix Optimization Problems. PhD thesis, National University of Singapore (2012)
17. Ding, C., Sun, D.F., Zhang, L.W.: Characterization of the robust isolated calmness for a class of conic programming problems. SIAM J. Optim. **27**(1), 67–90 (2017)
18. Dontchev, A.L., Rockafellar, R.T.: Characterizations of Lipschitz stability in nonlinear programming. In Mathematical Programming With Data Perturbations, Marcel Dekker, New York, 65–82 (1997)
19. Dontchev, A.L., Rockafellar, R.T.: Implicit Functions and Solution Mappings. Springer, New York (2009)
20. Dorsch, D., Gómez, W., Shikhman, V.: Sufficient optimality conditions hold for almost all nonlinear semidefinite programs. Math. Program. **158**(1), 77–97 (2016)
21. Eckstein, J., Silva, P.J.S.: A practical relative error criterion for augmented Lagrangians. Math. Program. **141**(1), 319-348 (2013)
22. Fernández, D., Solodov, M.V.: Local convergence of exact and inexact augmented Lagrangian methods under the second-order sufficient optimality condition. SIAM J. Optim. **22**(2), 384-407 (2012)
23. Fortin, M., Glowinski, R.: Augmented Lagrangian Methods: Applications to Numerical Solutions of Boundary Value Problems North-Holland, Amsterdam (1983)
24. Gafni, E.M., Bertsekas, D.P.: Two-metric projection methods for constrained optimization. SIAM J. Control Optim. **22**(6), 936–964 (1984)
25. Golshtein, E.G., Tretyakov, N.V.: Modified Lagrangians and Monotone Maps in Optimization. Wiley, New York (1989)
26. Han, D.R., Sun, D.F., Zhang, L.W.: Linear rate convergence of the alternating direction method of multipliers for convex composite programming. To appear in Math. Oper. Res. **42** (2017)
27. Han, D.R., Sun, D.F., Zhang, L.W.: Linear rate convergence of the alternating direction method of multipliers for convex composite quadratic and semi-definite programming. arXiv:1508.02134 (2015)
28. Hestenes, M.R.: Multiplier and gradient methods. J. Optim. Theory Appl. **4**(5), 303-320 (1969)
29. Higham, N.J.: Computing the nearest correlation matrix – a problem from finance. IMA J. Numer. Anal. **22**(3), 329–343 (2002)
30. Higham, N.J. and Strabić, N.: Bounds for the Distance to the Nearest Correlation Matrix. SIAM J. Matrix Anal. A. **37**(3), 1088–1102 (2016)
31. Ito, K., Kunisch, K.: The augmented Lagrangian method for equality and inequality constraints in Hilbert spaces. Math. Program. **46**(1), 341–360 (1990)
32. Izmailov, A.F., Kurennoy, A.S., Solodov, M.V.: A note on upper Lipschitz stability, error bounds, and critical multipliers for Lipschitz-continuous KKT systems. Math. Program. **142**(1), 591–604 (2013)
33. Jiang, K.F., Sun, D.F., Toh, K.-C.: Solving nuclear norm regularized and semidefinite matrix least squares problems with linear equality constraints. In *Discrete Geometry and Optimization*, Springer, 133–162 (2013)
34. Klatte, D.: Upper Lipschitz behavior of solutions to perturbed $C^{1,1}$ programs. Math. Program. **88**(2), 285–311 (2000)
35. Leventhal, D.: Metric subregularity and the proximal point method. J. Math. Anal. Appl. **360**(2), 681–688 (2009)
36. Li, X.D., Sun, D.F., Toh, K.-C.: QSDPNAL: A two-phase proximal augmented Lagrangian method for convex quadratic semidefinite programming. arXiv:1512.08872 (2015)

37. Luque, F.J.: Asymptotic convergence analysis of the proximal point algorithm. SIAM J. Optim. **22**(2), 277–293 (1984)
38. Mordukhovich,B.S., Sarabi, M. E.: Critical multipliers in variational systems via second-order generalized differentiation. Math. Program. doi:10.1007/s10107-017-1155-2 (2017)
39. Nilssen, T.K., Mannseth, T., Tai, X.-C.: Permeability estimation with the augmented Lagrangian method for a nonlinear diffusion equation. Computat. Geosci. **7**(1), 27–47 (2003)
40. Pataki, G., Tunçel, L.: On the generic properties of convex optimization problems in conic form. Math. Program. **89**(3), 449–457 (2001)
41. Pennanen, T.: Local convergence of the proximal point algorithm and multiplier methods without monotonicity. Math. Oper. Res. **27**(1), 170–191 (2002)
42. Powell, M.J.D.: A method for nonlinear constraints in minimization problems. In: Fletcher, R. (ed.) Optimization, 283-298. Academic, New York (1972)
43. Qi, H.D., Sun, D.F.: A quadratically convergent Newton method for computing the nearest correlation matrix. SIAM J. Matrix Anal. A. **28**(2) 360–385 (2006)
44. Robinson, S.M.: An implicit-function theorem for generalized variational inequalities. Technical Summary Report No. 1672, Mathematics Research Center, University of Wisconsin-Madison, 1976; available from National Technical Information Service under Accession No. ADA031952.
45. Robinson, S.M.: Some continuity properties of polyhedral multifunctions. Math. Program. Study. **14**, 206–214 (1981)
46. Robinson, S.M.: Constraint nondegeneracy in variational analysis. Math. Oper. Res. **28**(2), 201–232 (2003)
47. Rockafellar, R.T.: Convex Analysis. Princeton University Press, Princeton (1970)
48. Rockafellar, R.T.: Local boundedness of nonlinear monotone operators. Michigan Math. J. **16**(4), 397-407 (1969)
49. Rockafellar, R.T.: Conjugate Duality and Optimization. SIAM, Philadelphia (1974)
50. Rockafellar, R.T.: Augmented Lagrangians and applications of the proximal point algorithm in convex programming. Math. Oper. Res. **1**(2), 97–116 (1976)
51. Rockafellar, R.T.: Monotone operators and the proximal point algorithm. SIAM J. Control Optim. **14**(5), 877–898 (1976)
52. Rockafellar, R.T., Wets, R.J.-B.: Variational Analysis. Springer, New York (1998)
53. Shapiro, A.: First and second order analysis of nonlinear semidefinite programs. Math.Program. **77**(2), 301–320 (1997)
54. Shapiro, A.: Sensitivity analysis of generalized equations. J. Math. Sci. **115**(4), 2554–2565 (2003)
55. Shapiro, A., Sun, J.: Some properties of the augmented Lagrangian in cone constrained optimization. Math. Oper. Res. **29**(3), 479–491 (2004)
56. Sun, D.F., Sun, J., Zhang, L.W.: The rate of convergence of the augmented Lagrangian method for nonlinear semidefinite programming. Math. Program. **114**(2), 349–391 (2008)
57. Sun, J.: On Monotropic Piecewise Qudratic Programming. PhD thesis, University of Washington, Seattle (1986)
58. Toint, P.L.: Global convergence of a class of trust-region methods for nonconvex minimization in hilbert space. IMA J. Numer. Anal. **8**(2), 231–252 (1988)
59. Yang, L.Q., Sun, D.F., Toh, K.-C.: SDPNAL+: A majorized semismooth Newton-CG augmented Lagrangian method for semidefinite programming with nonnegative constraints. Math. Program. Comp. **7**(3), 1–36 (2015)
60. Zhao, X.Y., Sun, D.F., Toh, K.-C.: A Newton-CG augmented Lagrangian method for semidefinite programming. SIAM J. Optim. **20**(4), 1737–1765 (2010)
61. Zhou, Z.R., So, A.M.C.: A unified approach to error bounds for structured convex optimization problems. Math. Program. doi:10.1007/s10107-016-1100-9 (2017)

# Appendix.

## 1. Proof of Lemma 2.

*Proof* For notational simplicity, we write $s_{w,y} := -\mathcal{A}^* w - \mathcal{B}^* y - c$ for any $(w, y) \in \mathbb{W} \times \mathbb{Y}$. The "if" part is obviously true. To prove the "only if" part, we first show that there exist positive constants $\varepsilon$ and $\mu$ such that

$$-g^0(w, \bar{y}, s_{w,\bar{y}}) \geq -\sup(\mathrm{D}) + \mu \|w - \bar{w}\|^2, \quad \forall\, w \in \mathbb{B}_\varepsilon(\bar{w}). \tag{32}$$

It follows from $\bar{w} \in \operatorname{dom} h^*$ and the local strong convexity of $h^*$ that there exist positive constants $\varepsilon$ and $\mu$ such that

$$h^*(w) \geq h^*(\bar{w}) + \langle \nabla h^*(\bar{w}), w - \bar{w} \rangle + \mu \|w - \bar{w}\|^2, \quad \forall\, w \in \mathbb{B}_\varepsilon(\bar{w}) \cap \operatorname{dom} h^*.$$

By the convexity of $p^*$, we also get

$$p^*(s_{w,\bar{y}}) \geq p^*(s_{\bar{w},\bar{y}}) + \langle \varsigma, \mathcal{A}^*(-w + \bar{w}) \rangle, \quad \forall \varsigma \in \partial p^*(s_{\bar{w},\bar{y}}), \quad \forall w \in \mathbb{W}.$$

Hence, one can obtain the inequality (32) by adding the above two inequalities together and noting the KKT optimality condition (10). Now by shrinking $\varepsilon$ if necessary, we have, for any $(w, y, s) \in F_{D2} \cap \mathbb{B}_\varepsilon(\bar{w}, \bar{y}, s_{\bar{w},\bar{y}})$,

$$\begin{aligned}
-g^0(w, y, s) &\geq -g^0(y)/2 - g^0(w, \bar{y}, s_{w,\bar{y}})/2 \\
&\geq (-\sup(\mathrm{D}) + \kappa_2 \mathrm{dist}^2(y, \mathrm{SOL}_\mathrm{D}))/2 + (-\sup(\mathrm{D}) + \mu \|w - \bar{w}\|^2)/2 \\
&\geq -\sup(\mathrm{D}) + \min\{\kappa_2/2, \mu/2\}(\mathrm{dist}^2(y, \mathrm{SOL}_\mathrm{D}) + \|w - \bar{w}\|^2),
\end{aligned}$$

where in the second inequality the first term is due to the assumed quadratic growth condition for (D) at $\bar{y}$ with modulus $\kappa_2$, and the second term comes from (32). Finally, it follows that for any $\hat{y} \in \mathrm{SOL}_\mathrm{D}$ and any $(w, y, s) \in F_{D2} \cap \mathbb{B}_\varepsilon(\bar{w}, \bar{y}, s_{\bar{w},\bar{y}})$,

$$\begin{aligned}
\mathrm{dist}^2((w, y, s), \mathrm{SOL}_{D2}) &\leq \|w - \bar{w}\|^2 + \|y - \hat{y}\|^2 + \|s - s_{\bar{w},\hat{y}}\|^2 \\
&= \|w - \bar{w}\|^2 + \|y - \hat{y}\|^2 + \|\mathcal{A}^*(w - \bar{w}) + \mathcal{B}^*(y - \hat{y})\|^2 \\
&\leq (1 + 2\|\mathcal{A}^*\|^2)\|w - \bar{w}\|^2 + (1 + 2\|\mathcal{B}^*\|^2)\|y - \hat{y}\|^2,
\end{aligned}$$

which, with $\hat{y} := \Pi_{\mathrm{SOL}_\mathrm{D}}(y)$, implies

$$\mathrm{dist}^2((w, y, s), \mathrm{SOL}_{D2}) \leq (1 + 2\|\mathcal{A}^*\|^2)\|w - \bar{w}\|^2 + (1 + 2\|\mathcal{B}^*\|^2)\mathrm{dist}^2(y, \mathrm{SOL}_\mathrm{D}).$$

Thus, the quadratic growth condition (11) holds at $(\bar{w}, \bar{y}, s_{\bar{w},\bar{y}}) \in \mathrm{SOL}_{D2}$ for (D2).

**2. Proof of Proposition 1(c).**

*Proof* The convergence of $\{z^k\}$ under criterion (A) has been proven in Proposition 1. To establish the desired convergence rate, we first recall that the calmness of the mapping $T^{-1}$ at the origin for $z^\infty$ with modulus $\kappa$ asks for the existence of positive constants $\varepsilon$ and $\delta$ such that

$$\mathrm{dist}(z, T^{-1}(0)) \leq \kappa \|u\|, \quad \forall z \in T^{-1}(u) \cap \mathbb{B}_\delta(z^\infty), \quad \forall u \in \mathbb{B}_\varepsilon(0).$$

From parts (a) and (b) in Proposition 1 and the convergence of $\{z^k\}$, we obtain that

$$P_k(z^k) \in T^{-1}((z^k - P_k(z^k))/\sigma_k), \ \forall k \geq 0 \quad \text{and} \quad P_k(z^k) \to z^\infty \text{ as } k \to \infty,$$

which, imply the existence of a nonnegative integer $\bar{k}$ such that

$$\mathrm{dist}(P_k(z^k), T^{-1}(0)) \leq (\kappa/\sigma_k)\|z^k - P_k(z^k)\|, \quad \forall k \geq \bar{k}.$$

Now taking $\bar{z} = \Pi_{T^{-1}(0)}(z^k)$ in Proposition 1(b), we deduce that for any $k \geq 0$,

$$\begin{aligned}
\|z^k - P_k(z^k)\|^2 &\leq \|z^k - \Pi_{T^{-1}(0)}(z^k)\|^2 - \|P_k(z^k) - \Pi_{T^{-1}(0)}(z^k)\|^2 \\
&\leq \mathrm{dist}^2(z^k, T^{-1}(0)) - \mathrm{dist}^2(P_k(z^k), T^{-1}(0)).
\end{aligned}$$

Thus, it holds that

$$\mathrm{dist}(P_k(z^k), T^{-1}(0)) \leq \kappa/\sqrt{\kappa^2 + \sigma_{\bar{k}}^2}\,\mathrm{dist}(z^k, T^{-1}(0)), \quad \forall k \geq \bar{k}.$$

Hence, if criterion $(B)$ is also executed, we have that for any $k \geq \bar{k}$,

$$\begin{aligned}
&\|z^{k+1} - \Pi_{T^{-1}(0)}(P_k(z^k))\| \\
&\leq \|z^{k+1} - P_k(z^k)\| + \|P_k(z^k) - \Pi_{T^{-1}(0)}(P_k(z^k))\| \\
&\leq \eta_k \|z^{k+1} - z^k\| + \|P_k(z^k) - \Pi_{T^{-1}(0)}(P_k(z^k))\| \\
&\leq \eta_k(\|z^{k+1} - \Pi_{T^{-1}(0)}(P_k(z^k))\| + \|z^k - P_k(z^k)\|) + (\eta_k + 1)\|P_k(z^k) - \Pi_{T^{-1}(0)}(P_k(z^k))\| \\
&\leq \eta_k \|z^{k+1} - \Pi_{T^{-1}(0)}(P_k(z^k))\| + \left[\eta_k + (\eta_k + 1)\kappa/\sqrt{\kappa^2 + \sigma_k^2}\right] \operatorname{dist}(z^k, T^{-1}(0)).
\end{aligned}$$

Then the inequality in part (c) of Proposition 1 readily follows from the fact that $\operatorname{dist}(z^{k+1}, T^{-1}(0)) \leq \|z^{k+1} - \Pi_{T^{-1}(0)}(P_k(z^k))\|$ for any $k \geq 0$.

## 3. Proof of Proposition 3(b).

*Proof* By Lemma 1(b) and the convergence of $\{y^k\}$, we have $\operatorname{dist}(0, T_l(x^{k+1}, y^{k+1})) \to 0$ under criterion $(\widetilde{B})$. Therefore, by the upper Lipschitz continuity of $T_l^{-1}$ at the origin, we can derive, for $k$ sufficiently large,

$$\begin{aligned}
\operatorname{dist}((x^{k+1}, y^{k+1}), T_l^{-1}(0)) &\leq \kappa_l \operatorname{dist}(0, T_l(x^{k+1}, y^{k+1})) \\
&\leq \kappa_l (\operatorname{dist}^2(0, \partial f_k(x^{k+1})) + (1/\sigma_k^2)\|y^{k+1} - y^k\|^2)^{1/2} \\
&\leq (\kappa_l/\sigma_k)(1 + \eta_k'^2)\|y^{k+1} - y^k\|.
\end{aligned}$$

This completes the proof of this part.